\documentclass{report}
\usepackage[dvips]{graphics}
\usepackage{psfig}
\usepackage{amsfonts}
\usepackage{amssymb}
\usepackage{amsmath}
\usepackage{amsthm}
\usepackage{latexsym}
\newtheorem{thm}{Theorem}
\newtheorem{lem}[thm]{Lemma}
\newtheorem{prop}[thm]{Proposition}
\newtheorem{cor}[thm]{Corollary}
\theoremstyle{definition}
\newtheorem{defin}[thm]{Definition}
\newtheorem{conject}[thm]{Conjecture}
\newtheorem{remark}{Remark}

\def\Ad{\mathrm{Ad }}
\def\Ad_{\underline{\mathrm{Ad }}}
\def\Kb{K_P}
\def\kb{\mathfrak{k}_P}

\def\Mb{\mathrm{End}(V)_P}
\def\G2{\mathcal{G}^2_{con}}
\def\G4{\mathcal{G}^4_{con}}
\def\Poincare{Poincar\'{e} }
\def\C1{\mathcal{C}^1_{con}}
\def\C3{\mathcal{C}^3_{con}}
\def\k{\mathfrak{k}}
\def\grad{\mbox{grad }}
\def\Ck{\mathcal{C}^k_{con, A_0}}
\def\k{\mathfrak{k}}
\def\Gk{\mathcal{G}^{k+1}_{con}}
\def\Hk{\mathcal{H}^{k}_{con, A_0}}
\def\Hko{\mathcal{H}^{k}_{con, 0}}
\def\Lie{\mathrm{Lie}}
\def\Span{\mathrm{Span}}
\def\tr{\mathrm{tr}}

\begin{document}


\title{On the Holonomy of the Coulomb Connection over 3-manifolds with 
Boundary}
\author{William E. Gryc}
\date{August 2006}

\maketitle



\begin{abstract}
Narasimhan and Ramadas showed in \cite{NR} that the Gribov ambiguity was maximal for the product $SU(2)$ bundle over $S^3$. Specifically they showed that the holonomy group of the Coulomb connection is dense in the gauge group.  Instead of base manifold $S^3$, we consider here a base manifold with a boundary.  In this with-boundary case we must include boundary conditions on the connection forms.  We will use the so-called \emph{conductor boundary conditions} on connections.  With these boundary conditions, we will first show that the space of connections is a $C^\infty$ Hilbert principal bundle with respect to the associated conductor gauge group.  We will consider the holonomy of the Coulomb connection for this bundle.  If the base manifold is an open subset of $\mathbb{R}^3$ and we use the product principal bundle, we will show that the holonomy group is again a dense subset of the gauge group.
\end{abstract}




\tableofcontents

\setcounter{page}{1}
\pagenumbering{arabic}
\begin{chapter}{Introduction}
\begin{section}{Some Background to Yang-Mills Theory}
This thesis is concerned with classical Yang-Mills theory.  Our first task is to 
understand the very basics of this theory.  Physically, Yang-Mills theory 
models the strong nuclear force in the same way that Maxwell's equations model 
electromagnetic force.  To understand the mathematics, we will need to 
impose a lot of structure.  In this introduction, however, we will try 
to minimize the mathematical structure to get to the main idea, and refer 
the reader to Section 2 of \cite{BL} for the rigorous definitions.  Later we 
will carefully lay out the mathematics that we need.

Two basic concepts we will need are the Yang-Mills equation and the group of 
gauge transformations.

\begin{subsection}{The Yang-Mills Equation}

Consider the following situation:  Let $M$ be a compact 
oriented 3-dimensional manifold, and let $P\to M$ be a principal bundle over 
$M$ with compact connected structure group $K$.  Furthermore, assume that $K$ acts 
faithfully on a finite dimensional real (or complex) inner product space by 
isometries.  Thus we may view $K$ as a compact matrix subgroup of $U(V)$, and therefore the 
structure of $P$ induces a vector bundle $E:=P\times_K V \to M$ where 
$V=\mathbb{R}^n$ or $\mathbb{C}^n$ for some $n$.  

Given an Ehresmann connection $\omega_A$ on the principal bundle $P$, we obtain a Koszul 
connection $\nabla_A$ on the bundle $E$ whose curvature form 
$R_A$ is a vector-valued two-form on $M$.  The functional
\begin{equation}
\mathcal{YM}(\nabla_A) := \int_M |R_A|^2 dVol
\end{equation}
is well-defined for a certain norm.  Any connection that is a local minimum of 
the functional $\mathcal{YM}$ is called a \emph{Yang-Mills connection}.  One 
can show that $\nabla_A$ is a Yang-Mills connection if and only if $d^*_A R_A = 
0$. The equation 
\begin{equation*}
d_A^*R_A=0
\end{equation*} is called the \emph{Yang-Mills equation}.

The Yang-Mills equation has a 
similar structure to Maxwell's equations.  The Bianchi identity 
is that $d_AR_A =0$, so the Yang-Mills ``equations'' are
\begin{equation*}
d_AR_A =0,\mbox{    } d_A^*R_A=0.
\end{equation*} One can formulate Maxwell's equations in such a way that 
solving Maxwell's equations in the absence of a current is equivalent to 
finding a real-valued $2$-form $\eta$ on 
Minkowski space $\mathbb{R}^4$ such that
\begin{equation*}
d\eta = 0,\mbox{    } d^*\eta = 0.
\end{equation*} (For more on this construction, 
see the appendix to Chapter 2 and 10.2.8 in \cite{Dar}).  What makes solving 
the Yang-Mills equations difficult is the fact that the exterior derivative $d_A$ 
depends on the connection $\nabla_A$, and when written in local coordinates, this 
introduces a non-linear term of degree $3$.  In Maxwell's equations, the 
exterior derivative $d$ is independent of the form $\eta$, and the equations 
are linear.
\end{subsection}

\begin{subsection}{The Gauge Group}
There is a group $\mathcal{G}$, called the \emph{gauge group}, that acts on 
the set of connections and preserves the functional $\mathcal{YM}$.  So if $\nabla_A$ 
is a Yang-Mills connection and $g\in\mathcal{G}$, then $\nabla_A\cdot g$ is also a 
Yang-Mills connection.

The gauge group is fairly complicated to define in the general case, so here 
we will consider it locally.  Consider the principal bundle $P=\bar{O}\times K\to 
\bar{O}$, where $O$ is an open subset of $\mathbb{R}^3$ with smooth boundary, and 
$K$ is a compact matrix group.  Then we have the associated vector bundle $E=\bar{O}\times V\to \bar{O}$, where $V=
\mathbb{R}^n$ or $\mathbb{C}^n$, depending on whether the matrix group is 
real or complex.  A \emph{gauge transformation} $g$ is a mapping 
$g:\bar{O}\to K$ and the set $\mathcal{G}$ of gauge transformations is the 
\emph{gauge group}.  Given a section $\sigma$ of $E$ and $g\in\mathcal{G}$, 
we get a new section $g\cdot\sigma$ by
\begin{equation*}
(g\cdot\sigma)(x) = g(x)\sigma(x).
\end{equation*}  Given this left action on sections, we can get a right action 
on a connection $\nabla^A$ on the bundle $E$.  Indeed, we define the 
connection $\nabla^{A}\cdot g$ as
\begin{equation*}
(\nabla^{A}\cdot g)_X \sigma = g^{-1}\cdot\nabla^A_X(g\cdot\sigma)
\end{equation*}
for any vector $X\in T(O)$. It turns out that this action preserves connections that 
originally came from Ehresmann connections on the bundle $P$.
So if $\mathcal{C}$ is the space of connections on $P$, $\mathcal{G}$ acts on 
the right of $\mathcal{C}$.  Also, one can now check that if the connection 
$\nabla_A$ satisfies the Yang-Mills equation, then so does $\nabla_A\cdot g$.


Thus, a natural object to consider is the 
quotient space $\mathcal{C}/\mathcal{G}$, where $\mathcal{C}$ is the set of 
connections.  With some required modifications, the mapping $\mathcal{C}\to
\mathcal{C}/\mathcal{G}$ is an infinite dimensional principal bundle when $M$ 
is a compact 
$3$-manifold without boundary.  This bundle has been studied extensively in 
this case.  This thesis focuses on understanding this bundle when the 
underlying manifold $M$ has boundary.
\end{subsection}
\end{section}

\begin{section}{Conductor Boundary Conditions}
We will be concerned with boundary conditions on $p$-forms that have been dubbed 
\emph{conductor boundary conditions} by Gross in \cite{Gross}.  We say a form 
$\omega$ 
satisfies conductor boundary conditions if $i^*(\omega)=0$ where $i:\partial M
\to M$ is the inclusion map.  In other words $\omega(X_1\wedge\ldots\wedge X_p)
=0$ if $X_1,\ldots,X_p$ are all tangent to the boundary.  This is half of the 
``relative'' boundary conditions given by Ray and Singer in \cite{RS}, and the boundary conditions 
of the ``Dirichlet problem'' of Marini in \cite{Ma}.

We can extend the notion of conductor boundary conditions to connections.  
Given an Ehresmann connection $\omega_A$ on $P$, we consider the induced Koszul 
connection $\nabla_{A}$ on the vector bundle $E$.  Given a fixed connection 
$A_0$ on $P$, we say that a connection $\nabla_A$ satisfies {\it conductor boundary 
condtions with respect to $\nabla_{A_0}$} if the 1-form $\nabla_A - \nabla_{A_0}$ 
satisfies conductor boundary conditions.  Such a $\nabla_A$ equals 
$\nabla_{A_0}$ on the boundary in tangential directions, giving us a 
Dirichlet-like boundary condition (hence the terminology found in \cite{Ma}).
If we restrict our view to connections satisfying 
conductor boundary conditions with respect to a fixed $A_0$, we must change 
the gauge group so that it preserves the boundary conditions.  A gauge 
transformation $g\in\mathcal{G}$ is a section of certain bundle over $M$ with 
fibers diffeomorphic to $K$.  The \emph{conductor gauge group} $\mathcal{G}_
{con}$ consists of those gauge transformations $g$ such that $g|_{\partial M}
\equiv e$, where $e$ is the identity of $K$.  Note that the definition of 
$\mathcal{G}_{con}$ does not depend on the fixed connection $\nabla_{A_0}$.
\end{section}

\begin{section}{The Gribov Ambiguity and Holonomy of the Coulomb Connection}
The Gribov ambiguity comes up in the following setting, whose description is 
taken from \cite{Singer1}.  Physicists would like to 
compute a certain integral over $\mathcal{C}$ of the form
\begin{equation*}
\frac{\int_\mathcal{C} e^{-|R_A|^2}\{ \}\mathcal{D}A}{\int_\mathcal{C} e^{-|R_A|^2}\mathcal{D}A}.
\end{equation*}
This integral comes from the Feynman approach to quantum field theory. The problem with this integral is that the integrand in the numerator is constant on $\mathcal{G}$-orbits while the orbits are expected to have infinite 
measure. Morally, the integral should be taken over $\mathcal{C}/
\mathcal{G}$ and not $\mathcal{C}$.  So physicists do the following:  
instead of integrating over all of $\mathcal{C}$, take a continuous 
section $\sigma:\mathcal{C}/\mathcal{G}\to\mathcal{C}$ and integrate over 
$\sigma(\mathcal{C}/\mathcal{G})$, with an appropriate Jacobian 
weight factor from the change of variables.  They had a specific section in mind:
The infinite dimensional principal bundle $\mathcal{C}\to\mathcal{C}/
\mathcal{G}$ has its own connection called the \emph{Coulomb connection} 
with its horizontal subspaces given by
\begin{equation*}
H_A = \{\tau : \tau\mbox{ is a $\mathfrak{k}$-valued 1-form, } d_A^* \tau = 0\}.
\end{equation*}
We can define $\mathcal{S}_A=\{A+\tau: \tau\in H_A\}\subseteq\mathcal{C}$.  The physicists 
conjectured the following:

\begin{conject}\label{grib}
Fix a connection $\nabla_{A}\in\mathcal{C}$.  Then for every $\nabla_{A'}\in\mathcal{C}$, there exists a unique $g\in\mathcal{G}$ such 
that $\nabla_{A'}\cdot g\in\mathcal{S}_A$.
\end{conject}

We will see that locally this is true, meaning that given a small enough 
open set $\mathcal{O}\subseteq\mathcal{C}$ about $\nabla_{A}$, for every $\nabla_{A'}\in\mathcal{O}$ there exists a unique $g\in\mathcal{G}$ such that $\nabla_{A'}\in\mathcal{S}_A$.  However, Gribov showed in \cite{Gribov}\footnote{This is Singer's formulation of Gribov's result as found in \cite{Singer1}.}
\begin{thm}[The Gribov Ambiguity]
Suppose $M=S^4$, $P=S^4\times SU(2)\to S^4$, and the $\nabla_{A}$ is the flat connection (i.e. the Ehresmann connection whose 
horizontal subspaces are the tangent space to $S^4$ in $P$).  In particular, there exists a connection $\nabla_{A'}\ne \nabla_A$ in the $\mathcal{G}$-orbit of $\nabla_{A}$ such that $\nabla_{A'}\in\mathcal{S}_A$.  
\end{thm}
The ``ambiguity'' here is that given an $\nabla_{A}$, there might be multiple connections $\nabla_{A'}$ and $\nabla_{A''}$ that 
are gauge equivalent to $\nabla_{A}$ and are both in $\mathcal{S}_A$.  Hence, $\nabla_{A}$'s ``representative'' in 
$\mathcal{S}_A$ is ambiguous. 

Note that if Conjecture \ref{grib} were true, then the bundle $\mathcal{C}$ would be isomorphic to $\mathcal{G}\times 
\mathcal{S}_A$, and thus be a trivial bundle.  So, more generally, one can ask if if the bundle $\mathcal{C}$ is 
trivial, or equivalently ask if it allows any sections.  Both \cite{NR} and \cite{Singer1} show that
\begin{thm}[Generalized Gribov Ambiguity]
If $M=S^3$ or $S^4$, and $P=M\times SU(2)\to M$, then no continuous section $\sigma:\mathcal{C}/\mathcal{G}\to
\mathcal{C}$ exists.\footnote{Only \cite{Singer1} outlines the proof of the $M=S^4$ case.}
\end{thm}
The ``ambiguity'' here is that you cannot continuously choose a 
representative of each equivalence class of $\mathcal{C}/\mathcal{G}$.  
Hence, the physicists' idea of using a continuous section is mathematically impossible.

Narasimhan and Ramadas took the Gribov ambiguity a bit further for $M=S^3$ in the following sense.  
Given two connections $\nabla_A,\nabla_B\in\mathcal{C}$, they ask: how many points in the $\mathcal{G}$-orbit 
of $\nabla_B$ can be connected to $\nabla_A$ via horizontal paths with respect to the Coulomb connection?  The more points in the orbit can be connected by horizontal paths, the more ``ambiguous'' the Coulomb connection is.  The number of points that can be connected is the same as the number of elements in the holonomy group at $\nabla_A$ of the Coulomb connection, and thus this holonomy group becomes the main object of study.  Narasimhan and Ramadas show in \cite{NR} that if a 
certain metric is put on $S^3$ and the principal bundle considered is the product  
bundle $S^3\times SU(2)\to S^3$, then the holonomy 
group is dense in the connected component of the identity of $\mathcal{G}$. So they say that the ambiguity 
is maximal.  This also has some ramifications to physicists as 
is described in the introduction to \cite{NR}.

In this thesis, we address the question whether this maximal ambiguity 
holds if the base manifold is compact and \emph{with boundary} (unlike 
$S^3$) and $K$ is any compact semisimple matrix group.
\end{section}

\begin{section}{Summary of Results}
The aim of this thesis is to investigate the Gribov ambiguity when 
the base manifold $M$ is a compact 3-manifold with boundary, and the 
connections under consideration satisfy conductor boundary conditions 
with respect to a fixed connection $\nabla_{A_0}$.  We first need 
to prove that the corresponding bundle $\mathcal{C}\to\mathcal{C}/
\mathcal{G}$ is a $C^\infty$ principal bundle.  To this end, we will be 
using connection forms and gauge transformations of Sobolev classes $k$ and $k+1$ 
and denote them $\Ck$ and $\Gk$, respectively.  
Using standard techniques employed in \cite{AHS}, \cite{MV} and \cite{Par}, we will 
prove in Chapter \ref{bundle} that
\begin{prop} Suppose $M$ is a compact $3$-manifold with boundary, $P\to M$ 
is a principal bundle with a compact structure group $K$, and 
$k>3/2 + 1$.   Then $\Ck/\Gk$ is a $C^\infty$ 
Hilbert manifold, and $\Ck\to\Ck/\Gk$ is a principal bundle.
\end{prop}
With this proposition, it now makes sense to consider the holonomy 
group $\Hk$ of the Coulomb connection.  The infinite dimensional 
version of the Ambrose Singer 
theorem (see \cite{Magnot} for the statement of this theorem) 
tells us that $\Lie(\Hk)$ is generated by the image of the curvature 
form at certain points of the bundle.  Narasimhan and Ramadas use a similar but weaker fact 
to show that the holonomy group is dense in the connected component of the gauge group; indeed, they show that for a particular 
point $\omega\in\mathcal{C}$, the span of the image of the curvature of the Coulomb 
connection at $\omega$ is dense in the Lie algebra of the gauge group.  This 
leads directly to their result.  However, in Chapter 3 we prove that 
in our case the image of the curvature cannot be dense in the Lie algebra of the gauge group:
\begin{lem}  Let $M$ be a compact $3$-manifold with boundary and let $\nabla_{A}$ be a 
connection of Sobolev class $k$ for $k>3/2+1$ that satisfies conductor 
boundary conditions.  Define a set $\mathcal{L}_A$ as
\begin{equation*}
\mathcal{L}_A = \Span\{R_A(\alpha,\beta): \alpha,\beta\in H_A\}.
\end{equation*}
There exists a bounded nonzero operator $T_A:\Lie(\Gk)\to L^2(\kb|_{\partial M})$ such that 
$\mathcal{L}_A\subseteq \mathrm{ker}(T_A)$.
\end{lem}
Hence, the image of the curvature form at any \emph{fixed} point does not linearly generate the 
entire Lie algebra of $\Gk$.

We next specialize to the case where our principal bundle is the trivial bundle 
$P=\bar{O}\times K\to \bar{O}$, where $O\subset\mathbb{R}^3$ is an open subset with a smooth
boundary.  In this case, we can consider the flat connection.  Again, this is the Ehresmann connection 
whose horizontal subspaces are tangent to $\bar{O}$ in $P$.  We denote the corresponding 
Koszul connection on $\bar{O}\times V\to\bar{O}$ as $\nabla_0$ or simply $d$.  We then can show that 
the converse of the previous lemma holds for smooth functions if $\nabla_A = \nabla_0$.
\begin{lem}
Suppose we restricted the map $T_0$ above to smooth sections.  Then $\mathrm{ker}
(T_0|_{C^\infty})=\mathcal{L}_0\cap C^\infty$.
\end{lem}
We next consider the Lie algebra that $\mathcal{L}_0$ generates.
\begin{lem}
Let $g\in C^\infty(\partial M)$.  Then there exists $f\in [\mathcal{L}_0\cap C^\infty,
\mathcal{L}_0\cap C^\infty]$ such that $T_0(f)=g$.
\end{lem}
Using basic linear algebra and an argument of \cite{NR}, we prove
\begin{thm}
Let $f\in\Lie(\Gk)\cap C^\infty$.  Then $f$ is in the Lie algebra generated by $\mathcal{L}_0$.  
Hence, $\Hk$ is dense in the connected component of the identity of $\Gk$.
\end{thm}
Hence, in this special case, the maximal ambiguity of Narasimhan and Ramadas exists even when we are dealing with 
manifolds with boundary.
\end{section}
\end{chapter}
\begin{chapter}{The Bundle $\Ck / \Gk$}\label{bundle}
In this chapter we will prove that the mapping $\Ck\to\Ck/\Gk$ is 
indeed a $C^\infty$ vector bundle for $k>3/2 +1$.  This is a standard result, 
and versions of it have been proved in \cite{AHS}, \cite{MV}, \cite{NR}, and \cite{Par}.

%
\begin{section}{General Background and Notation}
$M$ will denote a 
compact oriented 3-dimensional Riemannian manifold with 
boundary, and $P\to M$ will denote a principal bundle with a semisimple compact 
structure group $K$. Furthermore, we assume that $K$ acts faithfully on a finite dimensional real 
(or complex) inner product space by isometries, and 
thus we view $K$ as a compact matrix group and a subgroup of $O(V)$ (or $U(V)$, respectively).  
Auxillary bundles also come into 
play.  The natural matrix multiplication of $K$ on $V:=\mathbb{R}^n$ (or V:= 
$\mathbb{C}^n$) induces a vector bundle $E:=P\times_{K} V$ (for the definition 
of these associated bundles, see Chapter 1.5 in \cite{KN}).  $K$ also acts on 
itself and $\mathfrak{k}$ via the adjoint representation, and thus we have the
corresponding bundles $\Kb := P\times_{K} K$ and $\kb := 
P\times_{K} \mathfrak{k}$.

Note that $\kb$ is a vector bundle, while $\Kb$ is not.  However, both $\kb$ 
and $\Kb$ are subbundles of the vector bundle $\Mb:=P\times_K 
\mathrm{End(V)}$, where again $K$ acts by the adjoint action.
 
Recall the exponential map 
$\exp :\mathfrak{k}\to K $.  Since $\Ad_ (k)\circ \exp = \exp\circ
\mathrm{Ad }(k)$, for any $k\in K$, 
we have an induced map $\exp :\kb\to\Kb$.

As $\mathrm{End}(V)$ acts on $V$ in an obvious way, fibers of $\Mb$ act 
on fibers of $E$.  Indeed, let $(p,T)_K\in\Mb$ and $(p,v)_K\in E$ be 
equivalence classes over the same point $x\in M$.  Then we define the 
action as
\begin{equation*}
(p,T)_K\cdot(p,v)_K=(p,Tv)_K.
\end{equation*}
It is easy to check that this is well-defined, 
and this induces a bundle isomorphism between $\Mb$ and 
$\mathrm{Hom}(E,E)$ over the identity.  
Viewing $\kb$ 
and $\Kb$ as subbundles of $\Mb$, fibers of $\kb$ and $\Kb$ also act on 
fibers of $E$. Similar reasoning also tells us that given two elements $\phi,
\psi$ in the same fiber of $\kb$, we can make sense of the Lie bracket 
$[\phi,\psi]$.  Indeed, if $\phi = (p,\phi')_K$ and $\psi = (p,\psi')_K$, 
then
\begin{equation*}
[(p,\phi')_K,(p,\psi')_K] = (p,[\phi',\psi'])_K
\end{equation*}
is well-defined by the Jacobi identity.

A Koszul connection $\nabla_A$ on $E$ induces a Koszul connection 
also called $\nabla_A^{Hom}$ on $\mathrm{Hom}(E,E)$ 
(see \cite{BL}, \cite{Dar} for more background on $\nabla_A^{Hom}$). Often, 
we will write $\nabla_A$ for $\nabla_A^{Hom}$ if it is clear we are using 
this induced connection.    
$\nabla_A^{Hom}$ induces a connection on $\kb$, and allows us to calculate 
$\nabla_A g$ for $g\in\Gamma(\Kb)$.  Note that $\nabla_A g$ is not 
necessarily a section of $\Kb$, but a section of $\mathrm{Hom}(E,E)$.  

Of special interest is the trivial bundle $P:=\bar{O}\times K\to\bar{O}$.  
In this case, the induced bundles $E$, $\kb$, and $\Kb$ are also direct 
products of the appropriate sort.  For example, $E=\bar{O}\times V\to \bar{O}$.  
In this case, the \emph{flat connection} on $P$ is the Ehresmann connection whose 
kernel is the tangent bundle of $\bar{O}$ embedded in the tangent bundle of $\bar{O}\times K$.  
Using parallel transport, one can check that the Koszul connection $\nabla_0$ induced on the product bundle $E$ is 
given by
\begin{equation*}
\nabla_0 \sigma = d \sigma\mbox{, $\sigma\in\Gamma(E)$}
\end{equation*}
where $d \sigma$ is the standard exterior derivative.  Thus we will often use $\nabla_0$ and $d$ interchangeably  
and call them the flat connection on $E$.

We are only concerned with connections $\nabla_A$ on $E$ that are 
induced by connections on $P$.  Such connections are called 
\emph{$K$-connections}, and one can show that $\nabla_A$ is a $K$-connection 
if and only if the local connection form is $\mathfrak{k}$-valued (see 
\cite{BL}).

As a vector bundle, we may equip $\kb$ with a metric.  Any Ad-invariant 
inner product on $\mathfrak{k}$ will induce a Riemannian metric on $\kb$.  In particular, 
we can use the trace inner product $(A,B)=\tr (A^*B)$ to induce a metric on 
$\kb$.  We now view $\kb$ 
as equipped with metric induced by the trace inner product on 
$\mathfrak{k}$.  Similarly, we equip the bundle $\Mb$ with the trace inner 
product.

For any vector bundle $\xi$ over $M$, we define 
the bundle $\Omega^j(\xi):=\mathrm{Hom}(\Lambda^j(M),\xi)$.  We call 
the elements of $\Omega^j(\xi)$ \emph{$\xi$-valued $j$-forms}, and 
generally call them \emph{vector-valued forms}.  (It would perhaps be better 
to call them ``vector-bundle valued forms,'' but this is not the standard 
terminology).  By convention, we have $\Omega^0(\xi):=\Gamma(\xi)$, the 
sections of $\xi$.  
Any connection on $\xi$ induces a connection on $\Omega^j(\xi)$ 
that involves the Levi-Civita connection on $M$.  See \cite{Dar} for more 
about these forms and this connection. 

There are certain operations we will like to define on forms.  Given any 
$\alpha\in\Omega^1(\kb)$ and $\phi\in\Gamma(\kb)$, we define the $1$-form 
$[\alpha,\phi]$ by
\begin{equation*}
[\alpha,\phi](X) = [\alpha(X),\phi],
\end{equation*}
for any $X\in TM$.  Also, given any $\alpha,\beta\in\Omega^1(\kb)$, we define the product 
$[\alpha\cdot\beta]\in \Gamma(\kb)$ in the following way:  Suppose 
locally $\alpha = \sum_i \alpha_i dx_i$, and $\beta = \sum_i \beta_i dx_i$, 
and the associated metric tensor is $\{g_{ij}\}$.  Then, we set
\begin{equation}\label{dotproduct}
[\alpha\cdot\beta] = \sum_{i,j} [\alpha_i,\beta_j]g^{ij},
\end{equation}
noting that $<dx_i,dx_j>=g^{ij}$ and the matrix $(g^{ij})$ is inverse to 
$(g_{ij})$ .  
One can verify that this globally defines $[\alpha\cdot\beta]$ as a section of $\kb$.

We will often be looking at the difference between two $K$-connections, 
and the following will be useful in looking at such differences. If 
$\nabla_{A_1}$ and $\nabla_{A_2}$ are 
$K$-connections, using the local characterization of $K$-connections, one 
can show that the difference 
$\nabla_{A_1}-\nabla_{A_2}$ is a $\kb$-valued $1$-form.  Furthermore, 
if we set $\alpha := \nabla_{A_1}-\nabla_{A_2}$, we have for any 
$\phi\in\Gamma(\kb)$
\begin{equation}\label{dprop}
(\nabla_{A_1}^{Hom}-\nabla_{A_2}^{Hom})(\phi)=[\alpha,\phi].
\end{equation}
Similarly, if $\beta\in\Omega^1(\kb)$, one can show that
\begin{equation}\label{d*prop}
((\nabla_{A_1}^{Hom})^*-(\nabla_{A_2}^{Hom})^*)(\beta) = -[\alpha\cdot\beta].
\end{equation}
The previous two equations are ubiquitous in what follows.  On sections we 
have $d_A=\nabla_A$ and on $1$-forms we have $d_A^*=(\nabla_A)^*$. We will 
use both notations interchangably on these respective domains.
The curvature $R_A$ of a $K$-connection $\nabla_A$ is a 
$\kb$-valued $2$-form.  All the facts asserted in this paragraph can be 
found in \cite{BL}.

Using \eqref{dprop}, we can show that a $K$-connection $\nabla_A^{Hom}$ is 
compatible with the metric on $\Mb$ induced by the trace inner product.  If we look locally, we have 
$\nabla_A^{Hom} = d + [A,\cdot]$, where $d$ is the flat connection with respect to a 
local coordinate system and $A$ is a local $\k$-valued $1$-form.  Then for a vector $X$ and local 
sections $S,T$ of $\Mb$, we have by the bilinearity of the inner product
\begin{equation*}
X\cdot<S,T> = <dS(X), T> + <S, dT(X)>.
\end{equation*}
Since we are locally working with matrices, we note that 
\begin{eqnarray*}
<[A(X),S],T> + <S,[A(X),T]> &=& \tr([A(X),S]^*T) +\tr(S^*[A(X),T]) \\
&=& \tr([S^*,A(X)^*] T) + \tr(S^*[A(X),T]) \\
&=& \tr(S^*A(X)^*T) - \tr(A(X)^*S^*T) + \\
& & \tr(S^*A(X)T) - \tr(S^*TA(X)) \\
&=& (-\tr(S^*A(X)T)) + \tr(A(X)S^*T)) + \\
& & \tr(S^*A(X)T) - \tr(S^*TA(X)) \\
&=& 0,
\end{eqnarray*}
where in the second to last line, we used the fact that $\mathfrak{k}\subseteq
\mathfrak{so}(V)$ (or $\mathfrak{su}(V)$).  Hence,
\begin{eqnarray*}
X\cdot<S,T> &=& <dS(X), T> + <S, dT(X)> \\
&=& <dS(X), T> + <S, dT(X)> + \\
& & <[A(X),S],T> + <S,[A(X),T]> \\
&=& <\nabla_A^{Hom} S, T> + <S,\nabla_A^{Hom} T>,
\end{eqnarray*}
proving metric compatibility.  Furthermore, a 
$K$-connection $\nabla_A$ on $E$ and the Levi-Civita connection on $M$ 
induce a connection on 
$\Omega^j(\kb)$ that is compatible with the induced metric on 
$\Omega^j(\kb)$.
\end{section}
%
%
%
%
\begin{section}{Sobolev Spaces of Connections and Gauge Groups}

We define Sobolev spaces of sections of vector bundles as Palais does in 
\cite{Palais}.  Using the notation of \cite{Palais}, the space $L^p_k(\xi)$ is the 
space of sections of $\xi$ with $k$ Sobolev derivatives under the $L^p$ norm, and $L^p_k(\xi)^0$ is the 
completion of $C^\infty_c(\xi|_{\mathrm{int}(M)})$ in the $L^p_k$ norm.  As usual, we define $H^k(\xi) := L^2_k(\xi)$, where the latter 
notation is what \cite{Palais} uses.  Also converting from Palais's notation, 
we put $H^k_0(\xi) := L^2_k(\xi)^0$.

Palais uses a local approach to define these Sobolev spaces.  However, it can be shown 
that given a smooth connection $\nabla_A$ on $\xi$, then the norm on $C^\infty$ sections $f$
\begin{equation*}
\|f\| = \|f\|_{L^p}+\sum_{i=1}^k \int_M <\nabla_A^i f,\nabla_A^i f>^{p/2} dVol
\end{equation*}
induces an equivalent norm on $L^2_k(\xi)$, and hence has the same completion.
We will use both
this global as well as Palais's local point of view of Sobolev spaces of sections.

Gross in \cite{Gross} has defined conductor boundary conditions on Sobolev 
spaces, which we will denote $H^k_{con}(\Omega^j(\xi))$ for appropriate vector 
bundles $\xi$, where $\Omega^j(\xi):=\mathrm{Hom}(\Lambda^j TM,\xi)$.\footnote{Marini in $\cite{Ma}$ 
has also defined these boundary conditions, although she calls them Dirichlet boundary conditions.} 
Specifically, 
we define the \emph{conductor Sobolev space $H^k_{con}(\Omega^j(\xi))$} for $k\ge 1$ as
\begin{equation}\label{defcon}
H^k_{con}(\Omega^j(\xi)):= \{\alpha\in H^k(\Omega^j(\xi)): \iota^*(\alpha)=0, \mbox{ where 
$\iota:\partial M\to M$ is the inclusion}\}.
\end{equation}
Since $k\ge 1$, $\alpha|_{\partial M}$ is defined in the trace sense, so $\iota^*(\alpha)$ is 
defined almost everywhere.  Again, $\iota^*(\alpha)=0$ is equivalent to saying that $\alpha$ vanishes 
on wedges of vectors $X_1\wedge\ldots\wedge X_j$, where all $X_i$ are tangent to $\partial M$.  
For a $0$-form $\sigma$ (i.e. a section $\sigma$), $\iota^*(\sigma)=0$ if and only if $\sigma|_{\partial M}=0$.
Hence, we see that
\begin{equation}
H^k_{con}(\xi) = H^1_0(\xi)\cap H^k(\xi)\mbox{, $k\ge 1$}\label{pB2}.
\end{equation}

In what follows we use $k>3/2 +1$ so we can use the multiplication theorem of 
Sobolev spaces (see Corollary 9.7 in \cite{Palais}).

Since we will be using conductor boundary conditions, we need a fixed smooth 
connection $\nabla_{A_0}$.  Set $\Ck := \nabla_{A_0} + H^k_{con}(\Omega^1(\kb))$. Note 
that all the connections in $\Ck$ will be equal to $\nabla_{A_0}$ in tangential 
directions on the boundary.  Also $\Ck$ is an affine space and thus seen 
to be a $C^\infty$-Hilbert manifold.  We will call any connection $\nabla_{A}$ 
\emph{$C^\infty$-smooth} if $\nabla_A-\nabla_{A_0}$ is a smooth section of $\kb$; 
in other words, $\nabla_{A}$ is a Koszul connection in the usual Riemannian geometry sense.  

\begin{prop}\label{mappings}
Suppose $k-1>3/2$, and $\nabla_{A}\in\Ck$.  Then for $1\le m\le k+1$, we have 
\begin{eqnarray*}
\nabla_A : H^{m+1}(\kb)\to H^m(\Omega^1(\kb))\\
\nabla_A^*: H^{m+1}(\Omega^1(\kb))\to H^m(\kb)
\end{eqnarray*}
are bounded linear transformations.
Also, $\nabla_A : H^{m+1}_{con}(\kb) \to H^m_{con}(\Omega^1(\kb))$ for $m=k-1$ and $m=k$. 
\end{prop}
\begin{proof}
Taking the global view of Sobolev spaces, we see that $\nabla_{A_0}$ maps $H^{m+1}(\kb)$ to 
$H^m(\Omega^1(\kb))$ and $\nabla_{A_0}^*$ maps $H^{m+1}(\Omega^1(\kb))$ to $H^m(\kb)$.  For $A\in\Ck$, if $1\le m\le k+1$, and $f\in H^{m+1}(\kb)$, we have
\begin{equation}\label{conductorpreserve}
\nabla_A f = \nabla_{A_0}f + [\nabla_A-\nabla_{A_0},f]\in H^m(\Omega^1(\kb))
\end{equation}
and that the mapping $f\to\nabla_A f$ is bounded by the multiplication theorem of Sobolev spaces.
Similarly, for $A\in\Ck$, and $1\le m\le k+1$, and $\alpha\in H^{m+1}(\Omega^1(\kb))$ we have 
\begin{equation*}
\nabla_A^* \alpha = \nabla^*_{A_0} \alpha - [\nabla_{A}-\nabla_{A_0}\cdot \alpha]\in H^m(\kb).
\end{equation*}
As for the last assertion of the proposition, note that if $f \in H^{m+1}_{con}(\xi)$, then $f\in C^1(\xi)$ and $f|_{\partial M}\equiv 0$.  Looking locally at a trivializing neighborhood at the boundary, we have 
$\nabla_A = d + A$, where $d$ is the flat connection with respect to the trivialization.  Since $A\in H^k_{con}(\Omega^1(\xi))$, we have $A\in C^1(\Omega^1(\xi))$.  Let $X$ be a tangential 
direction on the boundary.  Then since $f$ is constantly $0$ on $\partial M$, $df(X)=0$ on the boundary.  
Also, $[A(X),f]=0$ on the boundary since $f=0$ on the boundary.  Hence, globally, $\iota^*(\nabla_A f)\equiv 0$ on $\partial M$.   Hence by \eqref{defcon}, we have $\nabla_A f\in H^m_{con}(\Omega^j(\xi))$. 
\end{proof}

We now move onto the gauge transformations.  If $P$ were a trivial bundle, then the sections of $\Kb$ would be 
the gauge transformations described in the introduction.  Generally, the 
sections of $\Kb$ are the \emph{gauge transformations}.  The Sobolev regularity and boundary conditions we 
will need is set in the following definition:

\begin{defin}
Let $\nabla_{A_0}$ be a fixed smooth $K$-connection on $E$.  Let $g\in H^{k+1}(\Kb)$, with 
$g|_{\partial M}\equiv e$, where $e$ is the identity on $K$.  Then we say that $g\in\Gk$.
\end{defin}

The Sobolev space $H^{k+1}(\Kb)$ is defined as in \cite{MV} as the 
completion of smooth sections of $\Kb$ in the norm $H^{k+1}(\Mb)$.  This 
completion without 
the boundary conditions we will call $\mathcal{G}^{k+1}$, as it is called 
in \cite{MV}.

\begin{prop}
$\Gk$ is a closed topological subgroup of $\mathcal{G}^{k+1}$.
\end{prop}
\begin{proof}
Since the $H^{k+1}$ norm bounds the $L^\infty$ norm, we see that $\Gk$ is a closed subspace of $\mathcal{G}^{k+1}$.  
Since $e\cdot e= e$, and $e^{-1}=e$, we see that $\Gk$ is a closed topological subgroup of $\mathcal{G}^{k+1}$.
\end{proof}

$\Gk$ acts on $\Ck$ on the right in the 
following way.  Suppose that $\eta\in H^k_{con}(\Omega^1(\kb))$.  Then the 
action is
\begin{equation}
(\nabla_{A_0} + \eta)\cdot g = \nabla_{A_0} + g^{-1}\nabla^{Hom}_{A_0}g + 
\mathrm{Ad }(g^{-1})\eta.\label{gaction}
\end{equation}
Note that for $(\nabla_{A_0} + \eta)\cdot g$ to remain in $\Ck$, we need to have $g^{-1}\nabla^{Hom}_{A_0}g $ satisfy conductor boundary conditions.  The following proposition shows that this is the case.

\begin{prop}\label{gaugeEZ}
Suppose $g\in\Gk$.  Then we have $g^{-1}\nabla_{A_0}g\in H^k_{con}(\Omega^1(\kb))$, 
\end{prop} 
\begin{proof}
Let $g\in\Gk$.  Since $g\in H^{k+1}(\Kb)$, there exist smooth $g_n\in H^{k+1}(\Kb)$ such that $g_n\to g$ in 
$H^{k+1}(\Mb)$.  It is shown in \cite{MV} that inversion is continuous on $\mathcal{G}^{k+1}$.
Hence, $(g_n)^{-1}\to g^{-1}$ in $H^{k+1}(\Mb)$.  Since $\nabla_{A_0}$ is a smooth $K$-connection, we see that 
$g_n^{-1}\nabla_{A_0}^{Hom} g_n\in\Omega^1(\kb)$, and by the multiplication theorem,
\begin{eqnarray*}
||g^{-1}\nabla_{A_0}^{Hom} g - g_n^{-1}\nabla_{A_0}^{Hom} g_n||_{H^k} &\le& ||g^{-1}\nabla_{A_0}^{Hom} g - g_n^{-1}\nabla_{A_0}^{Hom} g||_{H^k}+\\ & &||g_n^{-1}\nabla_{A_0}^{Hom} g - g_n^{-1}\nabla_{A_0}^{Hom} g||_{H^k}\\
&\le& C(||g^{-1}-g_n^{-1}||_{H^{k+1}}||\nabla_{A_0}^{Hom}g||_{H^k} + \\
& & ||g_n^{-1}||_{H^{k+1}}||g_n - g||_{H^{k+1}})\to 0.
\end{eqnarray*}
Thus, $g^{-1}\nabla_{A_0}^{Hom} g\in H^k(\Omega^1(\kb))$.  

To show that $g^{-1}\nabla_{A_0}^{Hom} g\in H^k_{con}(\Omega^1(\kb))$ we proceed as in Proposition 
\ref{mappings}.  Locally, we have $\nabla_{A_0} = d + A_0$, where $A_0$ is a $C^\infty$-smooth $\mathfrak{k}$-valued 
$1$-form.  Let $X$ be a tangential direction.
Since $g\equiv e$ on $\partial M$, we have $dg(X)=0$. Also on the boundary, $[A_0(X),g]=[A_0(X),e]=0$, since 
$e$ commutes with everything.  Hence, globally, $\iota^*(\nabla_{A_0}^{Hom} g) \equiv 0$, and thus 
$\iota^*(g^{-1}\nabla_{A_0}^{Hom} g) \equiv 0$.  This proves that $g^{-1}\nabla_{A_0}^{Hom} g\in H^k_{con}(\xi)$, and 
thus $g\in\Gk$, as desired.
\end{proof}

We also note the following. If $\nabla_A\in \nabla_{A_0}+H^k(\Omega^1(\kb))$, 
$\eta:=\nabla_A-\nabla_{A_0}$, and $g\in\Gk$, then
\begin{equation*}
g^{-1}\nabla_A g = g^{-1}\nabla_{A_0} g + g^{-1}([\eta,g]) \in H^k_{con}
(\Omega^1(\kb))
\end{equation*}
since $[\eta,g]\equiv 0$ on the boundary ($g\equiv e$ on the boundary and thus commutes with 
everything). So $\nabla_{A_0}$ in Proposition \ref{gaugeEZ} can be replaced by 
any $H^k$ Sobolev connection.  In particular, it can be replaced by any smooth 
connection, and so $\Gk$ does not depend on the choice of the smooth 
$K$-connection $\nabla_{A_0}$.


Now we can state our desired result for this chapter:

\begin{thm}\label{bigpB}
Let $k>3/2+1$.  The quotient space $\Ck/\Gk$ is a $C^\infty$ Hilbert manifold, and $\pi: \Ck\to 
\Ck/\Gk$ is a principal bundle with structure group $\Gk$.
\end{thm}

The proof will be a straightforward adaptation of the work of Atiyah et al, 
Mitter et al, and Parker in \cite{AHS}, \cite{MV}, and \cite{Par}, 
respectively.  Before we can prove it, however, will need some basic tools 
which we will lay out in the next section.
  
\end{section}
%
%
%
%
\begin{section}{The \Poincare Inequality, Green Operators, and Other Necessities}

Before we start proving Theorem \ref{bigpB}, we need some facts about the 
Sobolev spaces we will be working with.  Many of these results will also be 
quite basic to the next chapter.  The most basic of these results is 
the following Sobolev-\Poincare inequality:

\begin{prop}\label{SPprop}
Let $\nabla_{A}$ be a smooth connection on a vector bundle $\xi\to M$ compatible with 
the metric.  Then there exists a $\kappa_p > 0$ such that for any 
$f\in L^p_1(\xi)^0$ with $1<p<\dim (M)$, we have
\begin{equation}
\|f\|_{L^p}\le\kappa_p\|\nabla_A f\|_{L^p}.\label{SP},
\end{equation}
where $\kappa_p$ does not depend on the connection, but does depend on $p$.
\end{prop}
\begin{proof}
First we consider real-valued functions on $M$, i.e. sections of the trivial 
bundle $M\times\mathbb{R}$.  We first want to prove that the only constant 
function 
in $L^p_1(M\times\mathbb{R})^0$ is the $0$ function.  So let $g$ be a constant 
function in $L^p_1(M\times\mathbb{R})^0$.  Using Theorem 9.3 in 
\cite{Palais}, we can extend the restriction map to a continuous map 
$R:L^p_1(M\times\mathbb{R})\to L^p(\partial M\times\mathbb{R})$.  Since 
$g\in L^p_1(M\times\mathbb{R})^0$, we have $g|_{\partial M}\equiv 0$.  But 
$g$ is constant, so $g\equiv 0$ everywhere, proving our assertion.

We can now use a standard argument in proving \Poincare inequalities.  
Variations can be found in the proof of Lemma 3.8 in \cite{He}, and 
in \cite{Evans}. We want to show that there exists a $\kappa_p > 0$ such that
 for any $g\in L^p_1(M\times\mathbb{R})^0$ with $1<p<\dim (M)$, we have
\begin{equation}
\|g\|_{L^p}\le\kappa_p\|dg\|_{L^p}.\label{HSP}
\end{equation}
Consider the set $\mathcal{H}:=\{g\in L^p_1(M\times\mathbb{R})^0: 
\|g\|_{L^p}=1\}$.   Note that we have just shown in the previous 
paragraph that no constant functions lie in $\mathcal{H}$.  Set
\begin{equation*}
C_p := \inf_{g\in\mathcal{H}} \|dg\|_{L^p}.
\end{equation*}
We will show that $C_p>0$, thereby proving \eqref{HSP} by setting $\kappa_p:=(C_p)^{-1}$.  Let 
$\{f_n\}\in\mathcal{H}$ attain the above infimum (in other words, 
$\lim_{n\to\infty}\|dg\|_{L^p}=C_p$).  Note that $\{f_n\}$ 
is bounded in $L^p_1$ norm and $L^p_1$ is reflexive (since $p>1$), 
so there exists a subsequence that we again call $\{f_n\}$ that 
has a weak limit $f$ in $L^p_1$.  By the Rellich-Kondrakov compactness 
theorem, since $1\le p < \dim{M}$, $\{f_n\}$ converges strongly in $L^p$ to this same 
function $f$.  The strong convergence 
shows that $f\in\mathcal{H}$, while the weak convergence shows that 
\begin{equation*}
\|df\|_{L^p}\le \liminf_{n\to\infty} \|df_n\|_{L^p}=C_p.
\end{equation*}
Since $f\in\mathcal{H}$, we know that $f$ is not constant, and thus $C_p=\|df\|_{L^p}$ 
is not $0$.  Hence, $C_p>0$, as desired, proving \eqref{HSP}.

Now let $f\in C^1_0(\xi|_{\mathring{M}})$.  Then the function 
$|f|$ is globally Lipschitz, so by Lemma 2.8 in \cite{He} we have 
$|f|\in L^p_1(M\times\mathbb{R})$.  Since $|f|$ is continuous and $0$ 
on $\partial M$, Theorem 5.5.2 in \cite{Evans} tells us that $|f|\in L^p_1(M
\times\mathbb{R})^0$.  Hence, \eqref{HSP} yields
\begin{eqnarray*}
\|f\|_{L^p}\le\kappa_p\|d|f|\|_{L^p}\le\kappa_p\|\nabla_A f\|_{L^p}.
\end{eqnarray*}
The second inequality is \emph{Kato's inequality}.  This inequality only requires that $\nabla_A$ is 
compatible with the metric.  For a proof of this inequality, see \cite{Par}.
Since $C^1_0(\xi|_{\mathring{M}})$ is dense in $L^p_1(\xi)^0$, the proposition 
has been proven.
\end{proof}

\begin{remark}
Note that while the connection $\nabla_{A}$ needs to be compatible with the metric 
to invoke Kato's inequality, 
the $\kappa_p$ above does not depend on the connection $\nabla_{A}$.  This allows us 
to let $\nabla_{A_\infty}$ be a $L^q$ connection for any $q\ge \dim(M)$ and still have 
(\ref{SP}) hold:  Let 
$\nabla_{A_n}\to \nabla_{A_\infty}$ in $L^q$ with $\nabla_{A_n}$ smooth.  Let $p^*$ be the Sobolev 
conjugate to $p$.  Invoking the boundary 
condition-free Sobolev inequality (found in \cite{He}), we have $f\in 
L^{p^*}$, and thus
\begin{eqnarray*}
\|f\|_{L^p} &\le &\kappa_p\|\nabla_{A_n}f\|_{L^p} \\
&\le &\kappa_p \|\nabla_{A_\infty}f\|_{L^p} + \|(\nabla_{A_n}-\nabla_{A_\infty})f\|_{L^p} \\
&\le &\kappa_p \|\nabla_{A_\infty}f\|_{L^p} + 
\|\nabla_{A_n}-\nabla_{A_\infty}\|_{L^{\dim(M)}}
\|f\|_{L^{p^*}} \\
&\to &\kappa_p \|\nabla_{A_\infty}f\|_{L^p},
\end{eqnarray*}
since $1/p=1/{p^*} +1/(\dim(M))$.  In particular, since $H^k(\kb)\subset C^1(\kb)\subset L^\infty(\kb)$, 
then \eqref{SP} holds for $\nabla_{A}\in\Ck$.
\end{remark}

The Sobolev-\Poincare inequality immediately tells us

\begin{cor}\label{gfree}
The action of $\Gk$ on $\Ck$ is free.
\end{cor}
\begin{proof}
Suppose $\nabla_{A_0} + \eta\in\nabla_{A_0} + H^k_{con}(\Omega(\kb))$, $g\in\Gk$, and 
$(\nabla_{A_0} +\eta)\cdot g = \nabla_A + \eta$.  By (\ref{gaction}), this means 
that
\begin{equation*}
(\nabla^{Hom}_{A_0} +\eta)g = \nabla^{Hom}_{A_0} g + [\eta,g] = 0
\end{equation*}
By (\ref{SP}), we have
\begin{equation*}
\|g-e\|_{L^2} \le \kappa_2\|(\nabla^{Hom}_{A_0} +\eta)g\|_{L^2} = 0
\end{equation*}
Hence, since $g$ is continuous, $g\equiv e$ and thus the corollary is proven.
\end{proof}

We are now in a position to show that $\Gk$ is a Hilbert Lie group.

\begin{prop}\label{liealg}
The exponential map takes $H^{k+1}_{con}(\kb)$ into $\Gk$ and is a local 
homeomorphism at $0$.
\end{prop}
\begin{proof}
In \cite{MV} it is shown that $\exp$ is a $C^\infty$ smooth map $\exp: H^{k+1}(\kb)\to\mathcal{G}^{k+1}$, 
as well as a local diffeomorphism, without boundary conditions.  Hence, for us to prove the proposition, 
we need only to show that $\exp$ maps $H^{k+1}_{con}(\kb)$ into $\Gk$, and for a neighborhood $U$ of the 
identity in $\mathcal{G}^{k+1}$, $\exp^{-1}=\log$ maps $U$ into $H^{k+1}_{con}(\kb)$.

To prove the first assertion, let $f\in H^{k+1}_{con}(\kb)$.  Since $k+1>3/2 + 2$ we have 
$f\in C^2(\kb)$ with $f|_{\partial M}\equiv 0$.  Then if $g:=\exp(f)$, we have that
$g\in C^2(\Kb)\subseteq C^2(\Mb)$, and 
$g|_{\partial M}\equiv e$, where $e$ is the identity element of $K$.  
Hence, by Proposition \ref{gaugeEZ}, $g\in\Gk$.

To prove the second assertion, since $\exp:H^{k+1}(\kb)\to\mathcal{G}^{k+1}$ 
is a local 
diffeomorphism between the spaces without boundary conditions, we need only 
to show that for small $\xi\in H^{k+1}(\kb)$, if $g:=\exp(\xi)\in\Gk$, 
then $\xi\in H^{k+1}_{con}(\kb)$.  Note that $\sup |\xi|\le C\|\xi\|_
{H^{k+1}}$. 
So for small enough $\|\xi\|_{H^{k+1}}$, we can use the fact that the ``pointwise'' 
map $\exp:\mathfrak{k}\to
K$ is local diffeomorphism at $0$ to say that since $g|_{\partial M}\equiv e$, 
we have $\xi|_{\partial M}\equiv 0$.  Since $\xi\in H^{k+1}(\kb)$, this 
vanishing on the boundary implies that $\xi\in H^1_0(\kb)\cap H^{k+1}(\kb)=
H^{k+1}_{con}(\kb)$.
\end{proof}
\begin{cor}\label{blgroup}
The group $\Gk$ is a Hilbert Lie group whose Lie algebra is identifiable with $H^{k+1}_{con}(\kb)$.  $\Gk$ acts 
smoothly on $\Ck$.
\end{cor}
\begin{proof}
The previous proposition shows that $\Gk$ is a local Hilbert Lie group, using logarithmic coordinates, with 
its Lie algebra identified with $H^{k+1}_{con}(\kb)$.  Exactly as in the proof of Theorem 2.18 of \cite{MV} (which shows 
$\mathcal{G}^{k+1}$ is a Hilbert Lie group), we can transport these coordinates throughout $\Gk$, and coordinate changes are smooth.  The details are in \cite{MV}. 

As for the smoothness of the action, it is shown in Proposition 3.12 of \cite{MV} that the action 
$\mathcal{C}^k\times\mathcal{G}^{k+1}\to\mathcal{C}^k$ is smooth.  The inclusions $\Ck\to\mathcal{C}^k$ and $\Gk\to\mathcal{G}^{k+1}$ are also smooth, so the composition $\Ck\times\Gk\to \mathcal{C}^k$ is smooth.  But we 
previously showed that action maps $\Ck\times\Gk$ into $\Ck$, as desired.
\end{proof}

Next on our agenda is to prove a version of Stokes' theorem.  In what follows, 
$M$ can be of any dimension.

\begin{lem}\label{Stokes'}
Let $M$ be a compact oriented Riemannian manifold with boundary, and 
$P\rightarrow M$ a principle $K$-bundle over $M$.  Let $u\in
\Omega^p(\kb)$ and $v\in\Omega^{p+1}(\kb)$. For any smooth connection $\nabla_{A}$, we 
have
\begin{eqnarray*}
\int_M (<d_A u, v> - <u,d^*_A v>) = \int_{\partial M} <\nu^*\wedge u, v> ,
\end{eqnarray*}
where $\nu$ is the outward pointing normal vector on the boundary, and $\nu^*$ 
is the Hilbert space adjoint of $\nu$ using the metric on $M$.
\end{lem}
\noindent In particular, if $u$ satifies conductor boundary conditions, then
\begin{eqnarray}\label{adjoint}
\int_M <d_A u, v> = \int_M <u,d^*_A v>.
\end{eqnarray}
To prove this we will prove a series of propositions.

\begin{prop}\label{frameind}
Let $M$ be a Riemannian manifold with dimension $n$, and let $\{e_i\}$ be an 
orthonormal frame field for $U\subseteq M$.  Let $E\rightarrow M$ be a vector 
bundle over $M$. Then for $u\in
\Omega^p(E)$ and $v\in\Omega^{p+1}(E)$,
\begin{eqnarray*}
X|_U = \sum_{i=1}^{n} <e_i^*\wedge u, v>e_i 
\end{eqnarray*}
defines a global vector field on $M$.
\end{prop}

\begin{proof}
Suppose $\{e_j\}$ and $\{f_j\}$ are two local orthonormal frame fields over
the same open set.  Then define the matrices $A = \{a_{ij}\}$ and $B=
\{b_{ij}\}$ by
\begin{eqnarray*}
e_j = \sum_i a_{ij} f_i ,\mbox{      }f_j = \sum_i b_{ij} e_i
\end{eqnarray*}
Then $B^{-1} = A$.  Also, note that since $e_j^*(f_i) = <e_j,f_i> = b_{ji}$, we
have
\begin{eqnarray*}
e_j^* = \sum_{i}b_{ji}f_i^* .
\end{eqnarray*}
Thus, since $B^{-1} = A$,
\begin{eqnarray*}
\sum_{j} <e_j^*\wedge u, v> &=& \sum_{i,j,k} <b_{ji}f_i^*\wedge u,v>a_{kj}f_k\\
 &=& \sum_{i,k}<f_i^*\wedge u,v>f_k(\sum_{j}a_{kj}b_{ji}) \\
 &=& \sum_{i,k}<f_i^*\wedge u,v>f_k\delta_{ik} \\
 &=& \sum_{i}<f_i^*\wedge u, v>f_i
\end{eqnarray*}
Thus, the definition of $X$ is frame independent, proving that X is a global vector field.
\end{proof}

An obvious approach to proving Lemma \ref{Stokes'} is to use a the regular 
Stokes' theorem, or a variation of it.  The variation we will use is the 
divergence 
theorem for manifolds.  This motivates the following proposition:

\begin{prop}\label{diveq}
If $X$ is the vector field defined by
\begin{eqnarray*}
X|_U = \sum_k <e_k^*\wedge u,v>e_k ,
\end{eqnarray*}
then
\begin{eqnarray*}
<d_A u, v> - <u,d^*_A v> = \mathrm{div}(X),
\end{eqnarray*}
where $\mathrm{div}(X)$ is taken with respect to the oriented volume form.
\end{prop}

\begin{proof}
We will modify Gross' proof in \cite{Gross} to our manifold case.
Choose a point $x\in M$ and let $\{e_j\}$ be a geodesic frame at that 
point.  First suppose $u$ and $v$ are locally of the form $u = \phi e_{j_1}^*
\wedge\ldots\wedge e_{j_p}^*$ and $v = \psi e_{k_1}^*\wedge\ldots\wedge 
e_{k_{p+1}}^*$, where $\phi$ and $\psi$ are local sections of our vector bundle
$\kb$.  Let $D$ be the Levi-Civita connection on $M$.  Then note at $x$ we 
have
\begin{eqnarray*}
D_{Y}e_{k_1}^*\wedge\ldots\wedge e_{k_{p+1}}^*(e_{l_1}\wedge\ldots\wedge
e_{l_{p+1}}) &=& Y(e_{k_1}^*\wedge\ldots\wedge e_{k_{p+1}}^*(e_{l_1}
\wedge\ldots\wedge e_{l_{p+1}})) \\
& &-\sum_{\alpha} e_{k_1}^*\wedge\ldots\wedge e_{k_{p+1}}^*(e_{l_1}\wedge \\
& &\ldots\wedge D_{Y}e_{l_\alpha}\wedge\ldots\wedge e_{l_{p+1}}) \\
&=& 0 - \sum_{\alpha}0 = 0.
\end{eqnarray*}
Thus, we have at $x$
\begin{eqnarray*}
d^*_Av &=& - \sum_{j}\iota_{e_j}
((\nabla_A)_{e_j}(\psi e_{k_1}^*\wedge\ldots\wedge e_{k_{p+1}}^*)) \\
&=&-\sum_{j}\iota_{e_j}(((\nabla_A)_{e_j}\psi) e_{k_1}^*\wedge\ldots\wedge 
e_{k_{p+1}}^* + \psi D_{e_j} e_{k_1}^*\wedge\ldots\wedge e_{k_{p+1}}^*) \\
&=&-\sum_{j}((\nabla_A)_{e_j}\psi) \iota_{e_j}(e_{k_1}^*\wedge\ldots\wedge 
e_{k_{p+1}}^*).
\end{eqnarray*}
Now we aim to obtain a similar expression for $d_A(u)$.  We first note that
at $x$ 
\begin{eqnarray*}
d(e_k^*)(e_i\wedge e_j) &=& e_i(e_k^*(e_j)) - e_j(e_k^*(e_i)) - 
e_k^*([e_i,e_j]) \\
& = & 0 - 0 - e_k^*(D_{e_i}e_j - D_{e_j}e_i) \\
& = & 0.
\end{eqnarray*}
Thus, at $x$
\begin{eqnarray*}
d(e_{j_1}^*\wedge\ldots\wedge e_{j_p}^*) = 0.
\end{eqnarray*}
So, at $x$
\begin{eqnarray*}
d_A u &=& (\nabla_A\phi)\wedge e_{j_1}^*\wedge\ldots\wedge e_{j_p}^* +
\phi d(e_{j_1}^*\wedge\ldots\wedge e_{j_p}^*) \\
&=& \sum_j ((\nabla_A)_{e_j}\phi)e_j^*\wedge e_{j_1}^*\wedge\ldots\wedge 
e_{j_p}^*.
\end{eqnarray*}
Let's define $\mu = e_j^*\wedge e_{j_1}^*\wedge\ldots\wedge e_{j_p}^*$ and 
$\omega = e_{k_1}^*\wedge\ldots\wedge e_{k_{p+1}}^*$.
Then using our expression for $d_Au$ and $d^*_Av$, we can write at $x$
\begin{eqnarray*}
<d_A u, v> - <u,d^*_A v> &=& \sum_j (<((\nabla_A)_{e_j}\phi)e_j^*\wedge \mu,v> 
+\\ 
& &-<u,((\nabla_A)_{e_j}\psi)\iota_{e_j}(\omega)>).
\end{eqnarray*}
Note that since $e_j^*\wedge\mu$ and $\omega$ are part of the orthonormal 
frame $\{e_{j_1}^*\wedge\ldots\wedge e_{j_{p+1}}^*\}$ for the bundle 
$\Omega^{p+1}(M)$, we have
\begin{eqnarray*}
<f e_j^*\wedge\mu,g \omega>&=&\sum <(f e_j^*\wedge\mu)(e_{l_1}\wedge\ldots
\wedge e_{l_{p+1}}),g \omega (e_{l_1}\wedge\ldots
\wedge e_{l_{p+1}})>_{\kb} \\
&=& <f,g>_{\kb} <e_j^*\wedge\mu,\omega>_{\Omega^{p+1}(M)}.
\end{eqnarray*}
We use the above and metric compatibility to get at $x$
\begin{eqnarray*}
<d_A u, v> - <u,d^*_A v> &=& \sum_j (<(\nabla_A)_{e_j}\phi,\psi> + 
<\phi,(\nabla_A)_{e_j}\psi>)\cdot\\
& &<e_j^*\wedge\mu,\omega> \\
&=& \sum_j e_j(<\phi,\psi>)<e_j^*\wedge\mu,\omega> \\
&=& \sum_j e_j(<\phi,\psi><e_j^*\wedge\mu,\omega>) \\
&=& \sum_j e_j(<\phi e_j^*\wedge\mu,\psi\omega>) \\
&=& \sum_j e_j(<u,\iota_{e_j}(v)>),
\end{eqnarray*}
where we also used the fact that $<e_j^*\wedge\mu,\omega>$ is 
constant.  So we have proven the proposition for all $u$ and $v$ with the 
special form $u = \phi e_{j_1}^*\wedge\ldots\wedge e_{j_p}^*$ and $v = 
\psi e_{k_1}^*\wedge\ldots\wedge e_{k_{p+1}}^*$. However, if 
$u$ and $v$ are arbitrary, then locally $u = \sum_i u_i$ and $v = 
\sum_j v_j$, where $u_i$ and $v_j$ are of the above special form.  So, in 
this general case, we have at $x$
\begin{eqnarray*}
<d_A u, v> - <u,d^*_A v> &=& \sum_{i,j} <d_A u_i, v_j> - <u_i,d^*_A v_j>
\\ &=& \sum_k \left(\sum_{i,j} e_k(<u_i,\iota_{e_k}(v_j)>)\right)\\
&=& \sum_k e_k(<u,\iota_{e_k}(v)>).
\end{eqnarray*}
However, since $\{e_k\}$ is a geodesic frame, this last sum is exactly the 
divergence of $X$ at the point $x$ (see Exercise 3.8a in \cite{doCarmo}).  So, 
we have
\begin{eqnarray*}
(<d_A u, v> - <u,d^*_A v>)(x) = \mathrm{div}(X)(x).
\end{eqnarray*}
But $x\in M$ was arbitrary, so we have our proposition.
\end{proof}
The final ingredient for Lemma \ref{Stokes'} is the divergence theorem for 
oriented Riemannian manifolds which states that
\begin{equation*}
\int_M \mathrm{div}(X) = \int_{\partial M} <\nu,X>.
\end{equation*}
Now, we can finish up proving the lemma.

\begin{proof}[Proof of Lemma \ref{Stokes'}.]
Let $\{e_i\}$ be an orthonormal frame including the boundary so that 
$\nu = e_1$. Then, for the vector field $X$ defined in Propostion \ref{diveq}, 
on $\partial M$ we have
\begin{eqnarray*}
<\nu,X> = <e_1,\sum_i<e_i^*\wedge u,v>e_i> = <e_1^*\wedge u,v> = 
<\nu^*\wedge u,v>.
\end{eqnarray*}
Thus, we have
\begin{eqnarray*}
\int_M (<d_A u, v> - <u,d^*_A v>) &=& \int_M \mathrm{div}(X) \\
&=& \int_{\partial M} <\nu,X> \\
&=& \int_{\partial M} <\nu^*\wedge u, v>,
\end{eqnarray*} 
which proves our lemma.
\end{proof}
Note that, like the Sobolev-\Poincare inequality, we can replace the smooth 
connection $\nabla_{A}$ with a Sobolev connection $\nabla_{A}$ in, for example, $\Ck$, and the 
smooth forms $u$ and $v$ with $H^1$ forms.

It will be essential for all that follows to have Green operators associated 
with every connection $\nabla_{A}\in\Ck$.  More specifically,
given a $K$-connection $\nabla_{A}\in\Ck$, we can define the Laplacian 
$\Delta_A=d^*_A d_A:H^{m+1}_{con}(\kb)\to H^{m-1}(\kb)$ for $1\le m\le k$.  
The regularity is correct by the following argument:  Since $\nabla_{A_0}$ is 
a smooth connection, clearly $\Delta_{A_0}$ is a bounded map from $H^m$ into 
$H^{m-2}$.  Suppose $h=\nabla_A-\nabla_{A_0}\in H^k_{con}(\kb)$.  Then for 
$f\in H^{m+1}$,
\begin{equation}\label{localLap}
\Delta_{A} f = \Delta_{A_0}f + [d_{A_0}^*h,f] - [h\cdot[h,f]] - 2[h\cdot 
d_{A_0}f].
\end{equation}
So, we have (allowing $||\cdot||_i$ to denote the $H^i$ norm)
\begin{eqnarray}
||\Delta_A f||_{m-1} &\le& ||\Delta_{A_0}f||_{{m-1}}+||[d_{A_0}^*h,f]||_
{{m-1}} \\
& &+ ||[h\cdot[h,f]]||_{m-1} + 2||[h\cdot d_{A_0}f]||_{m-1} \\
&\le& ||\Delta_{A_0}f||_{{m-1}}+C(||d_{A_0}^*h||_{k-1}||f||_{{m-1}} \\
& &+ ||h||_k||[h,f]||_{m-1} + 2||h||_k||d_{A_0}f]||_{m-1}) \\
&\le& ||\Delta_{A_0}f||_{{m-1}}+C(||h||_{k}||f||_{{m-1}} \\
& &+ ||h||_k^2||f||_{m-1} + 2||h||_k||f||_{m}) < \infty,\label{aok}
\end{eqnarray}
where we used the fact that $H^{m-1}$ is a $H^{k-1}$ module, which is the case 
since $k\ge m$ and $k-1\ge 3/2$.  Thus, $\Delta_A$ is bounded from $H^{m+1}_
{con}(\kb)$ to $H^{m-1}(\kb)$. Furthermore, we have
\begin{prop}\label{Greenexist}
Let $\nabla_{A}\in\Ck$ for $k-1\ge3/2$, and suppose $1\le m\le k$.  Then the mapping 
$\Delta_A:H^{m+1}_{con}(\kb)\to H^{m-1}(\kb)$ is an isomorphism.  Furthermore, 
if $f\in H^2_{con}(\kb)$ and $\Delta_A f\in H^{m-1}(\kb)$, then $f\in 
H^{m+1}_{con}(\kb)$ and 
$\Delta_A f$
\begin{equation}
\|f\|_{H^{m+1}}\le C(\|\Delta_A f\|_{H^{m-1}} + \|f\|_{H^0}).
\end{equation}
\end{prop}
We set $G_A:=(\Delta_A)^{-1}$ and call it the \emph{Green operator}.  
\begin{proof}
This proof will be a modification of a couple of standard arguments.  Our main 
sources will be \cite{Evans} and \cite{NR}.  There is probably a much cleaner 
way to prove this, but the author could not come up with one.  
Be forewarned: the following is nasty, brutish, and long.

First we note that if $\Delta_A f=0$ for some $f\in H^{m+1}_{con}(\kb)$, then 
by Lemma \ref{Stokes'}
\begin{equation*}
0=(f,\Delta_A f)_{L^2}=(d_A f, d_A f)_{L^2}.
\end{equation*}
So $d_A f = 0$. By the remark following Proposition \ref{SPprop}, since $f\in H^{1}_{0}(\kb)$, we have $f=0$.  So $\ker \Delta_A = 0$.  Most of the remainder of this proof will be showing $\Delta_A$ is onto.

We can define an bilinear form $B:H^1_0(\kb)\times H^1_0(\kb)\to\mathbb{R}$ as 
\begin{equation*}
B(u,v):=\int_M <d_A u, d_A v>.
\end{equation*}
Setting $h=\nabla_A - \nabla_{A_0}$ (and again allowing $\|\cdot\|_i$ to denote the $H^i$-norm), we 
have
\begin{eqnarray*}
|B(u,v)| &\le& ||d_A u||_0||d_A v||_0 \le (||d_{A_0}u||_0 + ||[h,u]||_0)(||d_{A_0}v||_0 + ||[h,v]||_0) \\
&\le& C(||u||_1 + ||h||_k||u||_0)(||v||_1 + ||h||_k||v||_0) \\
&\le& C||h||_k(||u||_1+||v||_1).
\end{eqnarray*}
So $B$ is bounded.  By Proposition \ref{SPprop}, $B$ is also coercive.  Let $f\in H^0(\kb)$ be arbitrary.  Then 
by the Lax-Milgram theorem (see Theorem 1 in Section 6.2.1 in \cite{Evans}), there exists a unique $u\in H^1_0(\kb)$ such 
that $B(u,v)=\int_M <f,v>$ for all $v\in H^1_0$.  The idea is that ``$\Delta_A u = f$'' (and if one considers 
$\Delta_A u\in H^{-1}$, then the quotes can be removed).  We now want to show that our solution $u$ is also in $H^2$.  To do this we will follow the proofs of Theorems 1 and 3 in Section 6.3 of \cite{Evans} closely.

We will show that $u\in H^2(\kb)$ by showing that $u\in H^2(\kb|_{V_k})$ for a finite cover $\{V_k\}$ of $M$.  First, 
let's take an interior trivializing coordinate neighborhood $U$ of $M$.  Let $D$ denote the flat connection on $U$.  
Then we set $\alpha:=\nabla_A - D\in H^k_{con}(\Omega^1(\kb|_U))$.  Let $\alpha =\sum_i \alpha_i dx_i$.
We also introduce the difference quotient $D^h_j g$ as 
\begin{equation*}
D^h_j g(x) = \frac{u(x+he_j)-u(x)}{h}.
\end{equation*}
where $e_k$ is the $k^{th}$ Euclidean direction vector.  Difference quotients behave much like derivatives.  For example, they is a Leibniz-type rule:
\begin{equation}\label{diffparts-1}
D^h_j(vw) = v^hD^h_j w + (D^h_jv)w,
\end{equation}
where $v^h(x)=v(x+he_j)$.
They also satisfy the integral equality
\begin{equation}\label{diffparts0}
\int_U <f, D^{-h}_j g> dx = -\int_U <(D^h_j f),g> dx  
\end{equation}
for $g$ with compact support in $U$.  Here, $dx$ is the Lebesgue measure on $\mathbb{R}^n$, and we need Lebesgue measure for 
\eqref{diffparts0} since we need the fact that it is translation invariant.  If $dVol$ is the volume measure induced by the metric of $M$, then $dVol=a dx$, where $a:=\sqrt{\det{(g_{ij})}}$.  Since both $a$ and $1/a$ are bounded on $U$, integrability under these two measures is equivalent.

Take open sets $V,W$ such that $\bar{V}\subseteq W\subseteq \bar{W}\subseteq{U}$, and choose $\zeta:U\to [0,1]$ so that $\zeta|_V\equiv 1$ and 
$supp(\zeta)\subseteq W$.  Set $v:=-D^{-h}_j(\zeta^2D^h_j u)$.  Note that by our choice of $\zeta$, we can extend $v$ by 
$0$ and have $v\in H^1_0(\kb)$.  Hence, we have
\begin{equation}\label{Green1}
\int_U <d_A u,d_A v> dVol = \int_U <f, v> dVol.
\end{equation} 
Set 
\begin{eqnarray*}
A&:=&\int_U <Du + [\alpha,u], Dv> dVol = \int_U <d_A u, Dv> dVol \\
B&:=&\int_U <f,v>dVol
-\int_U <Du,[\alpha,v]>dVol-\int_U <[\alpha,u],[\alpha,v]>dVol  \\
&=& \int_U <f,v>dVol-\int_U <d_A u,[\alpha,v]>dVol.
\end{eqnarray*}
Then \eqref{Green1} is equivalent to saying $A=B$.  Now 
\begin{eqnarray*}
A &=& -\sum_{i,l}\int_U <u_{x_i}+[\alpha_i,u], D^{-h}_j((\zeta^2D^h_j u)_{x_l})>g^{il}a dx \\
&=& \sum_{i,l}\int_U <D^h_j(g^{il}a (u_{x_i}+[\alpha_i,u])), (\zeta^2D^h_j u)_{x_l}> dx \\
&=& \sum_{i,l}\int_U (<(g^{il}a)^h D^h_j(u_{x_i}), \zeta^2D^h_j (u_{x_l})>)dx + \\
& & \sum_{i,l}\int_U (<(g^{il}a)^h D^h_j([\alpha_i,u]), \zeta^2D^h_j (u_{x_l})> + \\
& & <(g^{il}a)^h D^h_j(u_{x_i}+[\alpha_i,u]), 2\zeta\zeta_{x_l}D^h_j u> + \\
& & <D^h_j(g^{il}a) (u_{x_i}+[\alpha_i,u]), \zeta^2D^h_j (u_{x_l})> + \\
& & <D^h_j(g^{il}a) (u_{x_i}+[\alpha_i,u]), 2\zeta\zeta_{x_l}D^h_j u>dx).
\end{eqnarray*} 
Let the first sum of the last line equal $A_1$, and the second sum equal $A_2$.  Since the metric matrix $\{g_{ij}\}$ 
is positive and continuous, there exists a constant $\tilde{\theta}$ such that 
\begin{eqnarray*}
A_1 &>& \frac{\tilde{\theta}}{2}\int_U a^h \zeta^2\sum_i|D^h(u_{x_i})|^2dx \\
&\ge& \frac{\theta}{2}\int_U\zeta^2\sum_i|D^h(u_{x_i})|^2dx, 
\end{eqnarray*}
where $\theta = \tilde{\theta}\cdot\inf a$.  As for $A_2$, we have by Cauchy-Schwartz,
\begin{eqnarray*}
|A_2|&\le& C\int_U \zeta (|D^h_j\sum_i[\alpha_i,u]|\cdot|D^h_j\sum_{i}u_{x_i}| + |D^h_j\sum_{i}u_{x_i}|\cdot|D^h_ju| +\\
& & |D^h_j\sum_i[\alpha_i,u]|\cdot|D^h_ju| +|\sum_i[\alpha_i,u]|\cdot|D^h_ju| + |\sum_i u_{x_i}|\cdot|D^h_j\sum_{i}u_{x_i}|+ \\
& & |\sum_i[\alpha_i,u]|\cdot|D^h_j\sum_{i}u_{x_i}| + |\sum_i u_{x_i}|\cdot|D^h_ju|+|\sum_i[\alpha_i,u]|\cdot|D^h_ju|)dx.
\end{eqnarray*}
By the Peter-Paul inequality, for any $\epsilon>0$, we have
\begin{eqnarray*}
|A_2|&\le&\epsilon\int_U\zeta|D^h_j\sum_{i}u_{x_i}|^2 dx+\\
& & \frac{4}{\epsilon}\int_W |\sum_i u_{x_i}|^2 + |D^h_ju|^2 + |\sum_i[\alpha_i,u]|^2+|D^h_j\sum_i[\alpha_i,u]|^2dx.
\end{eqnarray*}
As in Theorem 3(i) of Section 5.8.2 in \cite{Evans}, we have
\begin{equation*}
\int_W |D^h_ju|^2dx \le C\int_U \sum_i |u_{x_i}|^2dx.
\end{equation*}
Letting $\epsilon=\theta/2$, we have
\begin{equation*}
|A_2|\le\frac{\theta}{2}\int_U\zeta\sum_i|D^h_j u_{x_i}|^2 dx+ C\int_U \sum_i|u_{x_i}|^2 + |[\alpha_i,u]|^2 +|[\alpha_i,u]_{x_l}|^2 dx.
\end{equation*}
So in sum, we have 
\begin{eqnarray*}
A &=& A_1 + A_2 \ge A_1 - A_2 \\
&\ge& \frac{\theta}{2}\int_U\zeta^2\sum_i|D^h(u_{x_i})|^2dx - C\int_U \sum_i(|u_{x_i}|^2+ |[\alpha_i,u]|^2+
 |\sum_{l}[\alpha_i,u]_{x_l}|^2)dx.
\end{eqnarray*}
Now we turn to $B$. We first look at the $L^2$ norm of $v$:
\begin{eqnarray*}
\int_U |v|^2 dx &=& \int_W |D^{-h}_j (\zeta^2 D^h_j u)|^2 dx \\
&\le& C \int_U |\sum_i (\zeta^2 D^h_j u)_{x_i}|^2 dx \\
&\le& C \int_U \zeta|D^h_j u|^2 + \sum_i\zeta|D^h_ju_{x_i}|^2 dx \\
&\le& C \int_U \sum_i|u_{x_i}|^2 + \sum_i\zeta|D^h_ju_{x_i}|^2 dx.
\end{eqnarray*}
Using the above and $<\gamma,[\beta,f]>=<[\gamma\cdot\beta],f>$ for $1$-forms $\gamma, \beta$ and $0$-form $f$, we have 
\begin{eqnarray*}
|B|&\le&C\int_U (|f|+|[\alpha\cdot Du]|+|[\alpha\cdot[\alpha,u]]|)\cdot |v| dx\\
&\le& \epsilon \int_U|v|^2 dx + C/\epsilon\int_U |f|^2+|[\alpha\cdot Du]|^2+|[\alpha\cdot[\alpha,u]]|^2dx \\
&\le& C \epsilon\int_U\sum_i\zeta|D^h_ju_{x_i}|^2 dx + C/\epsilon\int_U |f|^2+|[\alpha\cdot Du]|^2+|[\alpha\cdot[\alpha,u]]|^2+\\
& & \sum_i|u_{x_i}|^2dx.
\end{eqnarray*}
So choosing $\epsilon = \theta/(4C)$, we combine the $A$ and $B$ inequalities to obtain
\begin{eqnarray*}
(\int_V \sum_i |D^h_ju_{x_i}|^2dx)^{1/2} &\le& (\int_W \zeta\sum_i\zeta|D^h_ju_{x_i}|^2dx)^{1/2}\\
&\le& C||f||_0 +||[\alpha\cdot Du]||_0+ ||[\alpha\cdot[\alpha,u]]||_0+\sum_i||u_{x_i}||_0 + \\
& &\sum_i(||u_{x_i}||_0 + ||[\alpha_i,u]||_0+\sum_{l}||[\alpha_i,u]_{x_l}||_0)\\
&\le& C||\alpha||_k(||f||_0 + ||u||_1),
\end{eqnarray*}
where we used the fact that $H^1$ is a $H^k$ module in the last line, and $||\cdot||_k=||\cdot||_{H^k(\kb|_V)}$.  
Note that 
\begin{equation*}
\alpha = \nabla_A - D = (\nabla_A - \nabla_{A_0}) + (\nabla_{A_0}-D),
\end{equation*}
where the latter connection is smooth.  So, $||\alpha||_{H^k(\kb|_V)}\le C + ||\nabla_A-\nabla_{A_0}||_{H^k(\kb)}$.
So by Theorem 5.8(ii) of Section 5.8.2 in 
\cite{Evans}, we have that $u\in H^2(\kb|V)$ and
\begin{equation}\label{diffparts2}
||u||_{H^2(\kb|_V)}\le C(1+||\nabla_A-\nabla_{A_0}||_k)(||f||_0+||u||_1, )
\end{equation}
where we have switched back to $||\cdot||_k=||\cdot||_{H^k(\kb)}$.  Hence, we can replace $V$ above with any open set $V\subseteq M$ such that $\bar{V}\subseteq\mathring{M}$ and have \eqref{diffparts2} hold.  So, now we can understand 
$\Delta_A u$ as an function defined a.e. on $M$, and may deduce from 
Lemma \ref{Stokes'} that $\Delta_A u=f$ a.e. and all second derivatives exist as a.e. defined functions on $M$ (see the remark after Theorem 1 in Section 6.3 in \cite{Evans}.)  This observation is important for the following boundary considerations.

Let $U$ be a neighborhood of the boundary such that $U$ is the open unit ball $B_0(1)$ intersected with the upper half plane $\{x_3\ge 0\}$, and $\{x_3=0\}\cap U$ is the boundary portion of $U$.  We define another cut-off function $\zeta:\mathbb{R}^3\to[0,1]$ such that $\zeta|_{B_0(1/2)}\equiv 1$ and $supp(\zeta)\subseteq B_0(3/4)$.  Set $V:=B_0(1/2)$.  
For $j\in\{1,...,n-1\}$, we define $v:=-D^{-h}_j(\zeta^2 D^h_j(u))$.  One can check that
\begin{equation*}
v(x)=-\frac{1}{h^2}(\zeta^2(x-he_j)(u(x)-u(x-he_j))-\zeta^2(x)(u(x+he_j)-u(x)))
\end{equation*}
for $x\in U$.  Careful inspection above reveals that since $u|_{\partial M}=0$ in the trace sense, and $\zeta$ vanishes 
near the boundary of the ball, after extending by $0$ we have $v$ in $H^1_0(\kb)$. So,
\begin{equation*}
\int_U <d_A u, d_A v> dVol = \int_U <f,v> dVol.
\end{equation*}
By analogous estimates as the interior case, we have
\begin{equation*}
\int_V \sum_i |D^h_ju_{x_i}|^2dx \le C(1+||\nabla_A-\nabla_{A_0}||_k)^2(||f||_0+||u||_1)^2,
\end{equation*}
and so $u_{x_jx_i}\in H^2(\kb|_V)$ for all $i,j$ such that $i+j< 2\cdot 3$.  So the we need only consider $u_{x_3x_3}$.
Again, set $\alpha:=\nabla_A-D$, where $D$ is the flat connection.  Let $\Delta=D^*D$.  Then, recalling \eqref{localLap} we have 
\begin{eqnarray*}
u_{x_3x_3} &=& \Delta u -u_{x_1x_1}-u_{x_2x_2} \\
&=& f - [D^*h,u] + [\alpha\cdot[\alpha,u]] + 2[\alpha\cdot Du]-u_{x_1x_1}-u_{x_2x_2},
\end{eqnarray*}   
where equality is a.e.  Since everything on the right hand side is in $L^2(\kb|_V)$, so is 
$u_{x_3x_3}$.  Hence, $u\in H^2(\kb|_V)$.

Since we can cover $M$ with finitely many interior and boundary neighborhoods, we see that we have $u\in H^2(\kb)$, as we desired.  In sum, we have shown that given any $f\in L^2(\kb)$, there exists a $u\in H^2_{con}(\kb)$ such that $\Delta_A u = f$.  Our next job is to show that if $f\in H^{m-1}(\kb)$, then $u\in H^{m+1}$ for $1\le m \le k$.  Here we can use previous results, making things much easier.

The following argument is analogous to the one put forth in Section 3 of \cite{NR}.  Again suppose our $f$ from above is actually in $H^{m-1}(\kb)$.  Take any trivializing neighborhood $U$ of $M$ (interior or including the boundary).  Using our previous notation, by Theorem 5 of Section 6.3 in \cite{Evans} and the inequality following \eqref{localLap}, we have
\begin{eqnarray*}
||u||_{H^{m+1}(\kb|_U)} &\le& C(||\Delta u||_{m-1}+||u||_0)\\
&\le& C(||\Delta_A u||_{{m-1}}+||\alpha||_{k}||u||_{{m-1}} \\
& &+ ||\alpha||_k^2||u||_{m-1} + 2||\alpha||_k||u||_{m} +||u||_0)\\
&\le& C(||\Delta_A u||_{{m-1}} + ||u||_m + ||u||_0).
\end{eqnarray*}
where again $||\cdot||_i =||\cdot||_{H^{i}(\kb|_U)}$.  Interpolating on Sobolev norms, we have for any constant $\epsilon>0$
\begin{equation*}
||u||_m \le \epsilon||u||_{m+1} + C(\epsilon)||u||_0.
\end{equation*} 
So, choosing an appropriately small $\epsilon$, we have
\begin{equation*}
||u||_{H^{m+1}(\kb|_U)} \le C(||\Delta_A u||_{{m-1}} + ||u||_0).
\end{equation*}
Summing over a finite cover of $U$'s, we see that $u\in H^{m+1}_{con}(\kb)$.  Hence, $\Delta_A:H^{m+1}_{con}(\kb)\to H^{m-1}(\kb)$ is onto.  So by the Open Mapping Theorem, $\Delta_A$ is an isomorphism.  (One can also use compactness of $H^{m+1}$ in $L^2$ to directly show that the inverse is bounded, but we've done enough hard work for this proposition).
\end{proof}
\end{section}
%
%
%
%
\begin{section}{Proof of Theorem \ref{bigpB}}
Armed with the Lie algebra for $\Gk$ and Green operators $G_A$ for all $\nabla_{A}\in 
\Ck$, we can now start proving Theorem \ref{bigpB}.  

In our proof of Theorem \ref{bigpB}, we need to have some sort of slice 
lemma.  Informally, a ``slice'' is a chunk of the $\Ck$ to which our base $\Ck/\Gk$ will 
be locally diffeomorphic.  
Our slices will be modelled on the \emph{horizontal subspaces} $H_A$ at each $K$-connection $\nabla_{A}$.  
The horizontal subspace $H_A$ is defined as 
\begin{equation*}
H_A = \{\eta\in H^k_{con}(\Omega^1(\kb)): d_A^*\eta = 0\}.
\end{equation*}
For our slice lemma we need to have the following: for sufficiently small
$\eta \in H^k_{con, A_0}(\kb)$, there exists a unique gauge transformation 
$g\in\Gk$ 
so that $((\nabla_A + \eta)\cdot g)-\nabla_A \in H_A$.  
We would like to employ the implicit function 
theorem to find such gauge transformations, but clearly $\Gk$ is not a Banach 
space.  Instead, we consider the gauge algebra $H^{k+1}_{con}(\kb)$ and use it 
for the implicit function theorem.

We first prove the so-called ``local completeness'' of our action.  This proof is from \cite{AHS}:

\begin{prop}\label{lc}
Suppose $\nabla_{A}\in\Ck$.  Then there exists 
$\epsilon>0$ so that if $0<\epsilon_1 \le \epsilon$ there exists an $\epsilon_2>0$ such that for 
$\eta\in H^k_{con}(\Omega^1(\kb))$ with 
$\|\eta\|_{H^k}<\epsilon_1$ there exists a unique $g\in\Gk$ such that
$\|g-e\|_{H^{k+1}}<\epsilon_2$ and 
$(\nabla_A +\eta)\cdot g-\nabla_A\in H_A$.  In other words, 
$d_A^*((\nabla_A +\eta)\cdot g -\nabla_A)=0$.  Furthermore, $\epsilon_2$ can 
be made arbitrarily small by making $\epsilon_1$ sufficiently small.
\end{prop}
\begin{proof}
Consider the map $F:H^{k+1}_{con}(\kb)\times H^k_{con}(\kb)\to H^{k-1}(\kb)$ 
given by
\begin{equation}
F(X,\eta) = d_A^*(\exp(-X)d_A\exp(X) + \mathrm{Ad }(\exp(-X))\eta).
\end{equation}
F is just a composition of linear and bilinear maps, and the exponential 
map.  Thus, it is a smooth map of Banach spaces.  Also, $F(0,0)=0$.  To 
calculate the first partial derivative of $F$ at $(0,0)$, note that 
$G(X):= F(X,0) = d_A^*(\exp(-X)d_A\exp(X))$.  With the Baker-Campbell-Dynkin-Hausdorff 
formula, one can show that 
\begin{equation*}
D(\exp)(X)\xi = (\sum_{k=0}^\infty \frac{\mathrm{ad}^k(X)\xi}{(k+1)!})\exp(X).
\end{equation*}
So $D(\exp)(0)\xi=\xi$.  Thus, by Chain Rule and Proposition 14 in \cite{Lang1}, we have
\begin{eqnarray}
DG(0)\xi &=& d_A^*((D\exp)(0)\xi\cdot d_A \exp(0) + 
\exp(-0)d_A(D(\exp)(0)\xi))\label{prodrule}\\
&=& d_A^*d_A\xi.
\end{eqnarray}
To apply the implicit function theorem, we need to show that 
$d_A^*d_A=\Delta_A:H^{k+1}_{con}(\kb)\to H^{k-1}(\kb)$ is an isomorphism.  
This is exactly the statement of Proposition \ref{Greenexist}.  So applying the implicit 
function theorem we have:  There exists a $C^\infty$ mapping $X:N_A\cap H^k_{con}(\kb)
\to N_0\cap H^{k+1}_{con}(\kb)$, where $N_A$ is a neighborhood of $0$ in 
$H^k_{con}(\kb)$ and $N_0$ is a neighborhood of $0$ in $H^{k+1}_{con}(\kb)$
such that $X(\eta)$ is the unique member of $N_0$ satisfying
\begin{equation*}
d_A^*(\exp(-X(\eta))d_A\exp(X(\eta)) + \mathrm{Ad }(\exp(-X(\eta)))
\eta)=0.
\end{equation*}
By Proposition \ref{liealg}, $\exp: H^{k+1}_{con}\to\Gk$ is a local 
diffeomorphism at $0$.  So, given a small enough $\eta$, setting 
$g:=\exp(X(\eta))$ gives us the statement of the proposition.  This may require 
a shrinking of the neighborhood $N_A$, but this is possible since $X$ is a continuous 
map.  This last part of the proposition also follows from this continuity.
\end{proof}

Now we want so-called ``local effectiveness'' of our action.  This proof is again found in 
\cite{AHS} and, more directly \cite{Par} but with simplifications:
\begin{prop}\label{ge}
Suppose $\nabla_{A}\in\Ck$.  Then there exists an $\delta > 0$ so that
if $\|\eta_1\|_{H^k},\|\eta_2\|_{H^k}<\delta$, $\eta_1,\eta_2\in H_A$ and 
$\eta_1\not= \eta_2$ there exists no nontrivial $g\in\Gk$ such that 
$(\nabla_A +\eta_1)\cdot g = (\nabla_A +\eta_2)$.
\end{prop}
\begin{proof}
Suppose $\eta_1,\eta_2\in H_A$ and $(\nabla_A +\eta_1)\cdot g = (\nabla_A +
\eta_2)$, and $\|\eta_i\|_{H^k}<\delta$, for some $\delta>0$ that is yet to 
be determined. Set $\nabla_1 = \nabla_A +\eta_1$.  The idea here is to show that 
$\|g-e\|_{H^{k+1}}$ is small if both $\|\eta_1\|_{H^k}$ and $\|\eta_2\|_{H^k}$ 
are 
small, allowing us to invoke the uniqueness statement of Propostion \ref{lc}.
We have 
\begin{eqnarray}
\eta_2 &= & g^{-1}\nabla_{A} g + \mathrm{Ad }(g^{-1})\eta_1 \\
&= & g^{-1}\nabla_1 g \label{pB6}. 
\end{eqnarray}
Our main goal will be to show that $\|g-e\|_{H^{k+1}}$ is controlled by 
$\|g^{-1}\nabla_1 g\|_{H^k}$.  First note that since $K$ is compact, we have a constant 
$\Omega' \ge 1$ so that 
\begin{equation*}
\sup_{k\in K} |k| \le \Omega' < \infty
\end{equation*}
where the norm is induced by the trace inner product.  Hence $\sup_{g\in\Gk} ||g||_{L^\infty} \le 
\Omega'\cdot\mathrm{Vol}{M}=:\Omega$, since the norm on $\Mb$ is induced by 
the trace inner product.  Note that since $\nabla_1\in H^k$, we can apply 
(\ref{SP}) to obtain
\begin{equation}
\|g-e\|_{L^2} \le \kappa_2 \|\nabla_1(g-e)\|_{L^2} \le
\kappa_2\Omega \|g^{-1}\nabla_1 g\|_{L^2}.\label{pB7}
\end{equation}
By Proposition \ref{Greenexist} (note we could replace $\kb$ with $\mathrm{End}(V)_P$ and 
nothing would change in the proof of Proposition \ref{Greenexist}), since $\nabla_1 \in H^k$, 
we have
\begin{eqnarray}\label{generalGreen}
\|g-e\|_{H^{k+1}} &\le& C(\|\nabla_1^*\nabla_1 g\|_{H^{k-1}} + \|g-e\|_{L^2}).
\end{eqnarray} 
At first glace, one might expect that the constant $C$ above depends on $\eta_1$.  However, 
this is not the case.  Indeed,  
if we assume that $\delta < 1$ and thus $\|\eta_1\|_{H^k}<1$, we can use the inequality of 
Propostion \ref{Greenexist} with the connection $\nabla_A$ to get
\begin{equation}\label{generalGreen2}
\|g-e\|_{H^{k+1}} \le C_A(\|\nabla_A^*\nabla_A (g-e)\|_{H^{k-1}} + \|g-e\|_{L^2})
\end{equation}
The inequality \eqref{aok} with $\nabla_{A_0}$ replaced by $\nabla_A + \eta_1$ tells us that
\begin{eqnarray*}
\|\nabla_A^*\nabla_A (g-e)\|_{H^{k-1}} &\le & \|\nabla_1^*\nabla_1 (g-e)\|_{H^{k-1}} + C(||\eta_1||_{H^k}||g-e||_{H^{k-1}} \\
& &+ ||\eta_1||_{H^k}^2||g-e||_{H^{k-1}} + 2||\eta_1||_{H^k}||g-e||_{H^k}) \\
&\le & \|\nabla_1^*\nabla_1 (g-e)\|_{H^{k-1}} + C ||g-e||_{H^k},
\end{eqnarray*}
where we used the fact that $\|\eta_1\|_{H^k}<1$ on the last line.  Plugging the above into 
\eqref{generalGreen2} yields
\begin{equation*}
\|g-e\|_{H^{k+1}} \le C_A(\|\nabla_1^*\nabla_1 (g-e)\|_{H^{k-1}} +C||g-e||_{H^k}+ \|g-e\|_{L^2}).
\end{equation*}
A standard Sobolev norm interpolation then yields \eqref{generalGreen} with a constant $C$ that 
depends only on the fact that $\|\eta_1\|_{H^k}<1$.  

From \eqref{generalGreen} we have
\begin{eqnarray}
\|g-e\|_{H^{k+1}} &\le& C(\|g\|_{H^{k-1}}\|g^{-1}\nabla_1^*\nabla_1 g\|_{H^{k-1}} + \|g-e\|_{L^2})\\
&\le& C((\|g-e\|_{H^{k-1}}+\|e\|_{H^{k-1}})\|g^{-1}\nabla_1^*\nabla_1 g\|_{H^{k-1}} + \\
& & \|g-e\|_{L^2})\\
&\le & C((\|g-e\|_{H^{k-1}}+1)\|g^{-1}\nabla_1^*\nabla_1 g\|_{H^{k-1}} + \\
& & \|g^{-1}\nabla_1 g\|_{L^2}).\label{contle}
\end{eqnarray}
Note that
\begin{eqnarray*}
\|g^{-1}\nabla_1^*\nabla_1 g\|_{H^{k-1}}&\le& \|\nabla_1^*(g^{-1}\nabla_1 g)\|_{H^{k-1}} + 
\|\nabla_1 g^{-1}\cdot\nabla_1 g\|_{H^{k-1}} \\
&\le & C(\|g^{-1}\nabla_1 g\|_{H^{k}} + \|(\nabla_1 g^{-1})g\|_{H^{k-1}}\|g^{-1}\nabla_1 g\|_{H^{k-1}}) \\
&\le & C(\|g^{-1}\nabla_1 g\|_{H^{k}} + \|g^{-1}\nabla_1 g\|_{H^{k}}^2).
\end{eqnarray*}
In the above, we used the fact that $0=\nabla_1(g^{-1}g)=
(\nabla_1 g^{-1})g + g^{-1}\nabla_1g$ and Proposition \ref{mappings}.  Also, \emph{a priori} the constant $C$ should 
depend on $\eta_1$.  However, using reasoning similar to that which we used to show \eqref{generalGreen} tells us that $C$ 
depends only on the fact that $\|\eta_1\|_{H^k}<1$.  Since we assumed that $\delta < 1$, 
then by \eqref{pB6} we have $\|g^{-1}\nabla_1 g\|_{H^k}< 1$.  So we can remove the $\|g^{-1}\nabla_1 g\|_{H^{k}}^2$ 
term above.
Using the above and interpolation of Sobolev norms for $\|g-e\|_{H^{k-1}}$, we continue our inequality of \eqref{contle} with any $\epsilon_1 > 0$:
\begin{eqnarray*}
\|g-e\|_{H^{k+1}} &\le& C(\|g-e\|_{H^{k-1}} + 1)\|g^{-1}\nabla_1 g\|_{H^{k}}\\
&\le & \epsilon_1 \|g-e\|_{H^{k+1}}\|g^{-1}\nabla_1 g\|_{H^{k}} + C(\epsilon_1)\|g^{-1}\nabla_1 g\|_{H^{k}}
\end{eqnarray*}
Let's now assume that $\delta < 1/2$, which implies $\|g^{-1}\nabla_1 g\|_{H^{k}} < 1/2$.  So taking $\epsilon_1=1$, the 
above can be written as
\begin{equation}
\|g-e\|_{H^{k+1}} \le 2C(1)\delta\label{contle2}.
\end{equation}
Thus, taking $1>\epsilon>0$ as in Proposition \ref{lc}, we can choose
\begin{equation*} 
\delta <\min(\epsilon/(2 C(1)), 1/2,\epsilon).
\end{equation*}  
Then by \eqref{contle2}, we have 
$\|g-e\|_{H^{k+1}} < \epsilon$, so the uniqueness statement of Proposition 
\ref{lc} applies.  Since $\eta_1,\eta_2\in H_A$, this uniqueness tells us that 
$g\equiv e$.
\end{proof}
Now we can show that $\Ck/\Gk$ is a Hilbert manifold and prove local 
triviality of the bundle $\Ck\to\Ck/\Gk$.  Most of the following is exactly from 
\cite{MV} including most notation.

\begin{proof}[Proof of Theorem \ref{bigpB}]
Our quotient space $\Ck/\Gk$ first needs a topology. We  
give it the quotient topology under the projection $\pi:\Ck\to\Ck/\Gk$.
Fix a connection $\nabla_A\in\Ck$ and consider the mapping $F:H^k_{con}\times \Gk\to H^k_{con}$ 
given by
\begin{equation}
F(\eta,g)=(\nabla_A + \eta)\cdot g - \nabla_A .
\end{equation}
Then $F$ is continuous and $F(0,e)= 0$.  Hence, for a ball $B_{\delta}(A):=\{\nabla_A +\eta:
\|\eta\|_{H_k}<\epsilon \}$, there exists $\tilde{\epsilon_1}>0$ and $\tilde{\epsilon_2}>0$ so that 
for $\eta\in H^k_{con}$ and $g\in\Gk$ such that $\|\eta\|_{H^k}<\tilde{\epsilon_1}$ and 
$\|g-e\|_{H^{k+1}} < \tilde{\epsilon_2}$, then $F(\eta,g)\in B_{\delta}(A)$.  Set $\delta > 0$ to the 
$\delta$ in Propostion \ref{ge}.  Set $\epsilon_1>0$ so that it is less than $\min(\delta,\tilde{\epsilon_1})$, the $\epsilon$ in Proposition \ref{lc}, and so that the 
corresponding $\epsilon_2$ in Propostion \ref{lc} is less than $\tilde{\epsilon_2}$.  Set 
$\pi:\Ck\to\Ck/\Gk$ and $\mathcal{Q}_A:=\pi(B_{\epsilon_1}(A))$.
Consider the restriction $\pi_A := \pi|_{\mathcal{S}_A}$, where
\begin{equation*}
\mathcal{S}_A := {\pi^{-1}(\mathcal{Q}_A)\cap B_{\delta}(A)\cap (\nabla_A+ H_A)}.
\end{equation*}
Clearly $\pi_A$ maps into $\mathcal{Q}_A$.  We now show that this mapping is onto.
Given an equivalence class $[\nabla_{A'}]\in\mathcal{Q}_A$, we can assume without loss of generality 
that $\nabla_{A'}\in B_{\epsilon_1}(A)$.  By Proposition \ref{lc}, there exists $g\in\Gk$ with $\|g-e\|_{H^{k+1}}
<\epsilon_2$ and $F(\nabla_{A'}-\nabla_A,g)\in H_A$. Since $\nabla_{A'}\in B_{\tilde{\epsilon_1}}
(A)$ and $\|g-e\|_{H^{k+1}}<\epsilon_2 < \tilde{\epsilon_2}$, we see that
\begin{equation*}
\nabla_{A'}\cdot g = \nabla_A + F(\nabla_{A'}-\nabla_A,g)\in B_{\delta}(A).
\end{equation*}
Hence, $\nabla_{A'}\cdot g\in\mathcal{S}_A$ and $\pi_A(\nabla_{A'}\cdot g)=[A']$.  
So $\pi_A$ maps onto $\mathcal{Q}_A$.  Suppose $\nabla_{A_1},\nabla_{A_2}$ are in the domain of 
$\pi_A$ and $\pi_A(\nabla_{A_1})=\pi_A(\nabla_{A_2})$.  Then there exists $g\in\Gk$ such that 
$\nabla_{A_1}\cdot g = \nabla_{A_2}$.  Since $\nabla_{A_1},\nabla_{A_2}\in B_{\delta}(\epsilon)$, 
we can apply Proposition \ref{ge} to conclude that $\nabla_{A_1}=\nabla_{A_2}$.  Hence, $\pi_A$ is injective. 
Since $\mathcal{Q}_A$ has the quotient topology, the bijectivity of $\pi_A$ shows that it is a homeomorphism.
We will call its inverse $\sigma_A:\mathcal{Q}_A\to \mathcal{S}_A$.

We get a Hilbert 
manifold chart $\phi_A:\mathcal{Q}_A\to (\mathcal{S}_A-\nabla_{A})\subseteq H_A$ 
given by $\phi_A([\nabla_{A'}])=\sigma_A([\nabla_{A'}])-\nabla_{A}$.  It is easy to see that 
$(\mathcal{S}_A-\nabla_{A})$ is an open subset of $H_A$.  The next step is to show that 
coordinate changes are smooth.  To this end, we define a map $g_A:\pi^{-1}
(\mathcal{Q}_A)\to\Gk$ as follows: $g_A(\nabla_{A'})$ is the 
unique element of $\Gk$ so that 
\begin{equation}
\nabla_{A'}\cdot g_A(\nabla_{A'})^{-1} = \sigma_A([\nabla_{A'}]).\label{coorchange}
\end{equation}
$g_A(\nabla_{A'})$ exists and is unique by Corollary \ref{gfree}, and Propositions 
\ref{lc} and \ref{ge}.  
If $[\nabla_{A'}]\in\mathcal{Q}_{A_1}\cap\mathcal{Q}_{A_2}$, then 
we compute from (\ref{coorchange}) that 
\begin{eqnarray}
\phi_{A_1}([\nabla_{A'}]) &=& \sigma_{A_1}([\nabla_{A'}])-\nabla_{A_1} = \sigma_{A_1}([\sigma_{A_2}([\nabla_{A'}])])-\nabla_{A_1} \\
&=& \sigma_{A_2}([\nabla_{A'}])\cdot g_{{A_1}}(\sigma_{A_2}([\nabla_{A'}]))^{-1}-\nabla_{A_1} \\
&=&(\nabla_{A_2}+\phi_{A_2}([\nabla_{A'}]))\cdot\\
& & g_{A_1}(\nabla_{A_1} + (\nabla_{A_2}-\nabla_{A_1} + \phi_{A_2}([\nabla_{A'}])))^{-1}-\nabla_{A_1}.\label{coorchange2}
\end{eqnarray}
Since we know that gauge transformations are smooth, we need only show that 
$g_A$ is smooth for all $\nabla_{A}\in\Ck$ to show that this coordinate change is 
smooth.  This will come in the proof of local triviality of the quotient $\pi \Ck\to\Ck/\Gk$, 
which follows.

Note that when we show that the coordinate change is smooth, we will have 
a smooth map $\phi_{A_2}\circ\phi_{A_1}^{-1}$ from an open subset of $H_{A_1}$ 
to an open subset of $H_{A_2}$.  The first derivative of $\phi_{A_2}\circ
\phi_{A_1}^{-1}$ will thus provide an isomorphism from $H_{A_1}$ to $H_{A_2}$.

We want to show that a certain map $\Phi_A:\mathcal{Q}_A\times\Gk\to
\pi^{-1}(\mathcal{Q}_A)$ is a smooth diffeomorphism.  This map is given 
by $\Phi_A([\nabla_{A'}],g)=\sigma_A([\nabla_{A'}])\cdot g$.  Since gauge transformations are 
smooth, we see that $\Phi_A$ is smooth.  Also, $\Phi_A$ is a bijection with 
the inverse $\Phi_A^{-1}(\nabla_{A'})=(\pi(A'),g_A(\nabla_{A'}))$.  So, if we can show that 
$\Phi_A^{-1}$ is smooth, then $g_A$ will also be smooth making our coordinate 
change map (\ref{coorchange2}) smooth.  To consider the smoothness of 
$\Phi_A^{-1}$ we will look at $\Phi_A$ under coordinates and show that 
$\Phi_A$ is a local diffeomorphism at all points.

We know that $\mathcal{Q}_A$ is diffeomorphic to an open neighborhood 
$\tilde{\mathcal{S}_A}:=\mathcal{S}_A -\nabla_A$ in $H_A$ (since we haven't shown that $\Ck/\Gk$ is a manifold 
yet, to be correct we should replace $\mathcal{Q}_A$ in the domain of 
$\Phi_A$ with $\tilde{\mathcal{S}_A}$, prove smoothness of the inverse which then 
gives us that $\Ck/\Gk$ is a manifold, and then 
replace $\tilde{\mathcal{S}_A}$ with $\mathcal{Q}_A$ to give local triviality.  To 
avoid this extra confusing layer, we sweep this detail under the rug.) Given a 
fixed $g\in\Gk$, we have a 
neighborhood $M_g$ of the form $M_g = \{\exp(\xi)\cdot g: \xi\in H^{k+1}_{con}
(\kb),\|\xi\|_{H^{k+1}}<\epsilon\}$.  The set $V_\epsilon(0)=\{\xi\in H^{k+1}_
{con}
(\kb),\|\xi\|_{H^{k+1}}<\epsilon\}$ then provides coordinates for $M_g$.  
Finally, $\pi^{-1}(\mathcal{Q}_A)$ has coordinates under the mapping $\nabla_{A'}
\mapsto \nabla_{A'}-\nabla_{A}$.
So, we can rewrite $\Phi_A:\tilde{\mathcal{S}_A}\times V_\epsilon(0)\to\pi^{-1}
(\mathcal{Q}_A)-\nabla_{A}\subseteq H^k_{con}(\kb)$ as
\begin{eqnarray*}
\Phi_A(\tau,\xi) &=& g^{-1}\exp(-\xi)\nabla^{Hom}_A(\exp(\xi)g) + \mathrm{Ad}
(g^{-1}\exp(-\xi))(\tau) \\
&=& g^{-1}\exp(-\xi)(\nabla^{Hom}_A(\exp(\xi))g + g^{-1}\exp(-\xi)\exp(\xi)
(\nabla^{Hom}_Ag) + \\
& & \mathrm{Ad}(g^{-1}\exp(-\xi))(\tau) \\
&=&  \mathrm{Ad}(g^{-1})(\exp(-\xi)(\nabla^{Hom}_A(\exp(\xi))+\mathrm{Ad}
(\exp(-\xi))(\tau)) + g^{-1}\nabla^{Hom}_A g.
\end{eqnarray*} 
To use the inverse function theorem, we want to show that $(\Phi_A)_*$ is 
invertible at all points. Fixing
a $g\in\Gk$, and using the coordinates of $M_g$, we need only to consider 
the invertibility of $(\Phi_A)_*(\tau,0)$.  Since $\Phi_A$ restricted to 
the first variable is affine, we have 
\begin{equation*}
(\Phi_A)_{*1 (\tau,0)}(\eta) = \mathrm{Ad}(g^{-1})(\eta).
\end{equation*}
By a calculation similar to (\ref{prodrule}), we have
\begin{eqnarray*}
(\Phi_A)_{*2 (\tau,0)}(h) &=& \mathrm{Ad}(g^{-1})(\nabla_A^{Hom}h + -h\tau + 
\tau h) \\
&=&\mathrm{Ad}(g^{-1})(\nabla_{A'}^{Hom}h),
\end{eqnarray*}
where $\nabla_{A'}=\nabla_A + \tau$.  Adding up our partial derivatives yields
\begin{equation}
(\Phi_A)_{* (\tau,0)}(\eta,h) = \mathrm{Ad}(g^{-1})(\eta + 
\nabla_{A'}^{Hom}h)\label{phideriv}.
\end{equation}
To show $(\Phi_A)_{* (\tau,0)}$ is an isomorphism, we will first show that it 
has trivial kernel, and then show it is onto.  Also, we will drop the ``Hom'' 
from $\nabla_A^{Hom}$.

Suppose $(\Phi_A)_{* (\tau,0)}(\eta,h)=0$.  Since $\nabla_A^*\eta=0$, we 
have from (\ref{phideriv})
\begin{eqnarray*}
\Delta_Ah + \nabla_A^*[\tau,h] &=& \nabla_A^*(\nabla_A h + [\tau,h]) \\
&=& \nabla_A^*(\nabla_{A'}h) = \nabla_A^*(\nabla_{A'}h + \eta)\\
&=& \mathrm{Ad}(g)(\Phi_A)_{* (\tau,0)}(\eta,h) = 0.
\end{eqnarray*}
Applying the Green operator $G_A$ to the above yields
\begin{equation*}
h +G_A\nabla_A^*[\tau,h] = 0.
\end{equation*}
By the boundedness of $G_A:H^{k-1}(\kb)\to H^{k+1}_{con}(\kb)$ and $\nabla_A^*:H^k
(\kb)\to H^{k-1}(\kb)$, we have
\begin{eqnarray*}
\|h\|_{H^{k+1}} &\le & C\|\nabla_A^*[\tau,h]\|_{H^{k-1}} \\
&\le & C\|[\tau,h]\|_{H^k} \\
&\le & C\|\tau\|_{H^k}\|h\|_{H^{k+1}}
\end{eqnarray*}
For small enough $\tau$, the above implies that $h=0$, which in turn 
implies $\eta = 0$.  Thus, for small enough $\tau$, $\ker(\Phi_A)_{* (\tau,0)}$
is $0$.

Now we can move onto surjectivity.  Define a map $P_{A'}:H^k_{con}(\Omega^1(\kb))
\to H_{A'}$ as $P_{A'}(\omega) = (1-\nabla_{A'}G_{A'}\nabla^*_{A'})
(\omega)$. We can rewrite (\ref{phideriv}) as
\begin{equation}
(\Phi_A)_{* (\tau,0)}(\eta,h) = \mathrm{Ad}(g^{-1})(\nabla_{A'}(h
+ G_{A'}\nabla^*_{A'}\eta) + P_{A'}\eta).\label{phideriv2}
\end{equation}
We have written $(\Phi_A)_{* (\tau,0)}$ in the form 
$(\Phi_A)_{* (\tau,0)}:\tilde{\mathcal{S}_A}\times V_\epsilon(0)\to
\mathrm{Ad}(g^{-1})(\mathrm{Im}(\nabla_{A'}) \oplus H_{A'})$.  It is easy 
to see that $\mathrm{Im}(\nabla_{A'}) \oplus H_{A'}$ is indeed a direct 
sum and that $\mathrm{Im}(\nabla_{A'}) \oplus H_{A'}=H^k_{con}(\Omega^1(\kb))$ (see, 
for example, \cite{Gross} and \cite{NR}). Note that since $g\in\Gk$, $\mathrm{Ad}(g^{-1})$ maps 
$H^k_{con}$ to itself isomorphically.  Hence, to prove surjectivity, we must 
show that for every $h_0\in H^{k+1}_{con}$ and $\eta_0\in H_{A'}$, we 
have a (unique) solution to
\begin{equation*}
(\Phi_A)_{* (\tau,0)}(\eta,h) = \mathrm{Ad}(g^{-1})(d_{A'}h_0 + \eta_0).
\end{equation*}
Consider the function $H:H^k_{con}(\kb)\oplus H_A \oplus\tilde{\mathcal{S}_A}\to 
H^k_{con}(\kb)$ given by
\begin{equation*}
H(\tau,\eta_0, \eta) = \eta - d_{A'}G_{A'}[\tau\cdot\eta]-\eta_0.
\end{equation*}
Note that $H(0,0,0)=0$, $H$ is continuous and linear in the last two 
variables, and $H_{^*3(0,0,0)}=-\mathrm{Identity}$.  One can also show that 
$H$ is $C^1$ (see \cite{Gross} for example). 
So the implicit 
function theorem says that there exists an $\epsilon>0$ so if $\|\tau\|_{H^k}, 
\|\eta_0\|_{H^k}<\epsilon$, there exists an 
$\eta(\tau,\eta_0)$ such that $H(\tau,\eta_0,\eta(\tau,\eta_0)) = 0$.
Let $\eta_0\in H_A$ be arbitrary, and $\tau\in H^k_{con}(\kb)$ so that 
$\|\tau\|_{H^k}<\epsilon$.  Choose $N>0$ so that 
$\|(1/N)\eta_0\|_{H^k}<\epsilon$. Then using linearity in the last two 
variables we have $H(\tau,\eta_0,N\eta(\tau,(1/N)\eta_0)) = 0$.
Since $\eta:=N\eta(\tau,(1/N)\eta_0)\in H_A$, we have
\begin{eqnarray*}
\eta_0 &= &\eta - d_{A'}G_{A'}d^*_{A'}\eta \\
&= & P_{A'}(\eta).
\end{eqnarray*}
So we have a solution for $\eta_0$.  
Now set $h$ to
\begin{equation*}
h := h_0 - G_{A'}d^*_{A'}\eta.
\end{equation*}
Since $h_0$ satisfies boundary conditions and $G_{A'}$ maps $H^{k-1}$ into 
$H^{k+1}_{con}(\kb)$, we have $h\in H^{k+1}_{con}(\kb)$.  Furthermore
\begin{eqnarray*}
d_{A'} h_0 & = & d_{A'}(h+G_{A'}d^*_{A'}\eta).
\end{eqnarray*}
Thus, by (\ref{phideriv2}), we have found a solution to $(\Phi_A)_{* (\tau,0)}
(\eta,h) = \mathrm{Ad}(g^{-1})(d_{A'}h_0 + \eta_0)$.  Thus we 
have surjectivity for small $\tau$.  Hence $\Phi_A$ is a local diffeomorphism at all points, 
and therefore a diffeomorphism.  Local triviality is thus proven, and we 
have finally shown that $\Ck/\Gk$ is a Hilbert manifold and $\Ck\to\Ck/\Gk$ is 
a principal bundle.
\end{proof}
\end{section}
\end{chapter}
\begin{chapter}{The Holonomy of the Coulomb Connection}
Now that we know that the bundle $\Ck\to \Ck/\Gk$ is a principal bundle, 
we can consider holonomy group.  The connection we will consider is the called 
the \emph{Coulomb connection} whose connection form at $\nabla_A$ is defined as $G_A d^*_A$.  Then 
the corresponding horizontal at $\nabla_{A}$ is $H_A$.  Recall the definition of $H_A$ as
\begin{equation*}
H_A = \{\alpha\in H^k_{con}(\Omega^1(\kb)): d_A^*\alpha = 0\}.
\end{equation*}
This connection is natural in the sense that $H_A$ is the $L^2$ orthogonal 
complement to the vertical vectors at $\nabla_{A}$.  Indeed, one can show that given $\gamma\in 
\Lie(\Gk)=H^{k+1}_{con}(\Omega^1(\kb))$, the fundamental vector field associated to $\gamma$ is  
$d_A \gamma$. Hence the vertical vectors are those vectors of the form $d_A\gamma$ for some 
$\gamma\in\Lie(\Gk)$ (see \cite{Gross} or \cite{NR}). 
By the same reasoning as the proof of Lemma 7.1 in \cite{NR}, the Coulomb connection is indeed 
a connection on $\Ck\to\Ck/\Gk$. 

We begin our investigation of the holonomy group by considering the image of 
the curvature form $\Omega$ of the Coulomb connection.  Let $\mathcal{R}_A$ be the 
curvature form of the Coulomb connection at $\nabla_A$.  By the same calculation 
in the proof of Lemma 7.2 in \cite{NR}, we have 
\begin{equation}
\mathcal{R}_A(\alpha,\beta)=-2G_A([\alpha\cdot\beta]),\mbox{ for } \alpha,\beta\in 
H_A.
\end{equation} 
In this investigation, certain types of coordinates at the boundary are useful, 
and are the subject of the next section.
%
%
%
%
\begin{section}{Coordinates at the Boundary}
Consider the following system of coordinates at the boundary that satisfy the following:
\begin{enumerate}
\item[A1.] $\partial/\partial x_n$ is orthogonal to $\partial/\partial x_1,\ldots,\partial/\partial x_{n-1}$ on the boundary.\label{A1}
\item[A2.] $\partial/\partial x_n$ has norm 1 everywhere.\label{A2}
\item[A3.] $\partial/\partial x_n$ is the inward pointing normal vector on the 
boundary.\label{A3}
\end{enumerate}
We describe such a coordinate system as \emph{Type A}. Fortunately, this 
defintion is not in vain, as such coordinates always exist:
\begin{prop}\label{RSC}
Let $M$ be a Riemannian $n$-manifold with boundary.  A coordinate system 
$\{x_1,\ldots,x_n\}$ satisfying A1-A3 above exists around 
each point of the boundary of $M$.
\end{prop}
\begin{proof}
The following construction is based on \cite{Milnor1} and \cite{RS}, and \cite{RS} uses this 
type of coordinates.  Let $p$ 
be a point on the boundary.  Take a chart on the set $U$ 
near $p$ with 
coordinates $(y_1,\ldots,y_n)$ and image in the upper half space so that $y_n^{-1}(0)\cap U=\partial M\cap 
U$.  Then the 
function $u(y_1,\ldots,y_n):=y_n$ satifies $u^{-1}(0)=\partial M\cap U$.  Let 
$\nu$ be the inward pointing normal.  For any point $p'\in\partial M$, 
$d_{p'}u(\nu)=c\cdot\partial/\partial y_n u = c > 0$ where $c=||\partial/
\partial y_n||^{-1}$.  Let $X$ be a local vector field 
that is dual to the $1$-form $du$, i.e. $X = \mbox{grad } u$.  Since $u$ has 
no critical points, $X$ never vanishes.  So we may set $Y=X/(||X||)$.  Then, 
as in \cite{Milnor1}, 
we may consider the flow of $Y$, denoted $\Phi$.  As in \cite{RS},
$\Phi:[0,\delta)\times \partial M\cap U\to U$ is a diffeomorphism for some 
$\delta >0$.  We now define coordinates via this diffeomorphism:  Let $x_i :
\partial M \to \mathbb{R}$ be coordinates for $\partial M\cap U$, with 
inverse $\psi$.  
Define $x_i:U\to\mathbb{R}$ be defined as 
\begin{equation*}
x_i(\Phi(t,q))=x_i(q).
\end{equation*}
and define $x_n:U\to\mathbb{R}$ as
\begin{equation*}
x_n(\Phi(t,q))=t.
\end{equation*}
These $\{x_1,\ldots,x_n\}$ are coordinates on $U$.  Now, note that
\begin{eqnarray*}
\frac{d}{dt}(u\circ\Phi_t(q))|_{t=t_0} &=& du(X/{||X||})(\Phi_{t_0}(q)) \\
&=& (X,X/{||X||})(\Phi_{t_0}(q))=||X||(\Phi_{t_0}(q)).
\end{eqnarray*}
Hence, since $u\circ\Phi_0(q)=u(q)=0$, we have
\begin{equation*}
u(\Phi_t(q)) = \int_0^t\|X\|(\Phi_s(q))ds.
\end{equation*}
Using the above, we see that on the boundary,
\begin{eqnarray*}
du\left(\frac{\partial}{\partial x_i}\right) = 0 = dx_n\left(\frac{\partial}
{\partial x_i}\right), \mbox {for $i<n$} \\
du\left(\frac{\partial}{\partial x_n}\right) = ||X|| = ||X||dx_n\left(\frac
{\partial}{\partial x_n}\right).
\end{eqnarray*}
Hence, 
\begin{equation}\label{okgo}
du|_{\partial M}=||X||\cdot dx_n|_{\partial M}.
\end{equation}
Now we can start showing that our properties are satisfied.  For a function 
$f$ on $U$
\begin{eqnarray*}
\frac{\partial}{\partial x_n} f = \frac{d}{dt} f\circ\Phi(t,\psi(x_1,\ldots,
x_{n-1})) = \frac{X}{||X||}\cdot f.
\end{eqnarray*}
Thus, $\partial/\partial x_n = X/||X||$.  In particular, $||\partial/\partial 
x_n ||=1$, satisfying property A2 above.  For $i<n$, on the 
boundary we have by (\ref{okgo})
\begin{eqnarray*}
0&=&\frac{d}{dt}(x_i(q)) = \frac{d}{dt}(x_i(\Phi(t,q))) = \frac{1}{||X||}
dx_i(\grad u)\\
&=& \frac{1}{||X||}du(\grad x_i) \\
&=& dx_n(\grad x_i) = (\grad x_n, \grad x_i) \\
&=& g^{in}. 
\end{eqnarray*}
Since $\{g^{ij}\}$ is a symmetric matrix, this implies that $g_{in}=0$ on the 
boundary, proving property A1 above.  As for property A3, since $g_{in}=0$, 
$\frac{\partial}{\partial x_n}$ is normal to the boundary.  So it is 
either inward or outward pointing.  Since $(x_1,\ldots, x_n)$ is a chart on the upper half plane, we have 
that $\frac{\partial}{\partial x_n}$ is inward pointing by definition, completing the proof.
\end{proof}

We will also have the occasion to use a slightly different type of 
coordinates on the boundary.  If in the proof of Proposition \ref{RSC} we 
instead let $Y=X/||\grad X||^2$ where $X = \grad u$, and let $\Psi$ be 
the corresponding flow, then 
\begin{equation*}
\frac{d}{dt}(u\circ \Psi(t,q))=du(X/||X||^2) = (X,X/||X||) = 1.
\end{equation*}
So $u(\Psi(t,q))=t$.  So we can let $y_n=u$ and $y_i(\Psi(t,q)):=y_i(q)$ $i=1,\ldots
n-1$, where $(y_1,\ldots,y_{n-1})$ is a chart on the boundary.  Then following reasoning 
similar to the proof of Proposition \ref{RSC} gives us coordinates that 
satisfy
\begin{enumerate}
\item[B1.] $\partial /\partial y_n$ is orthogonal to $\partial /\partial y_1,\ldots\partial /\partial y_{n-1}$ 
everywhere.
\item[B2.] $\partial /\partial y_n$ is a positive (perhaps nonconstant) multiple of the 
inward pointing normal.
\end{enumerate}
The last condition follows from the fact that $du(\nu) > 0$ on the boundary.  
We creatively describe such coordinates as \emph{Type B}.  In what follows, 
we will use $\{x_1,x_2,x_3\}$ to denote Type A coordinates and $g_{ij}$ to 
denote the associated metric tensor. For Type B coordinates, 
we use $\{y_1,y_2,y_3\}$ and $\{h_{ij}\}$.

Also, if $h_{ij}$ is the metric tensor of a Type B coordinate system, note 
that by condition B1 we have
\begin{equation*}
h_{33}=\frac{1}{h^{33}}=\frac{1}{||\grad u||^2}.
\end{equation*}

\end{section}
\begin{section}{Mean Curvature}
It turns out that the mean curvature of the boundary comes into play in our 
characterization of the image of the curvature form.  We give a brief 
explanation of mean curvature and calculate it using Type A and Type B coordinates.
We use \cite{doCarmo} for this background.  Let $f:S\to M$ be an immersed 
submanifold, and let $\nabla$ be the Levi-Civita connection on $M$.
The \emph{second fundamental form} $B$ is a mapping 
$B:T_pS\times T_pS\to (T_pS)^\perp$ given by
\begin{equation}
B(x,y):= (\nabla_x Y)^N,
\end{equation}
where $Y$ is any local extension of $y$, and $Z^N$ is the normal component of 
a vector $Z\in T_pM$ with respect to $S$.  While not immediately apparant, it 
one can verify that $B$ is well-defined, symmetric, and a bilinear mapping 
of $C^\infty(S)$ modules (see \cite{doCarmo}).  Given a fixed $\eta\in (T_pS)^\perp$, the Riesz 
representation theorem gives us a mapping $S_\eta :T_pS\to T_pS$ satisfying 
\begin{equation}
(S_\eta(x),y) = (B(x,y),\eta).
\end{equation}
One can show
\begin{equation}\label{3}
S_\eta(x)=-(\nabla_x \tilde{\eta})^T,
\end{equation}
where $\tilde{\eta}$ is a local extension of $\eta$, and $Z^T$ is the tangent 
component of a vector $Z\in T_pM$.

The trace of this operator $S_\eta$ is 
important.  $S$ is called \emph{minimal} if $\tr(S_\eta)=0$ for all $\eta\in
(T_pS)^\perp$ and $p\in S$.  If $S$ is an oriented hypersurface and 
$\eta\in(T_pS)^\perp$ has norm $1$ and is pointing in the direction 
corresponding to the orientation, then 
\begin{equation}
H:=\frac{1}{\dim(S)}\tr(S_\eta)
\end{equation}
is called the \emph{mean curvature} of $f$. For us, the relevant immersion is 
$\iota:\partial M\to M$ and the normal vector will be the outward pointing 
normal which we denote $-\nu$ (so the inward pointing normal is still $\nu$).

A certain quantity will come up often when working with Type A coordinates on 
our $3$-manifold $M$. 
Let $\{g_{ij}\}$ be the metric tensor in a Type A coordinate system, and let $a=
\sqrt{\det(g_{ij})}$.  Note that $a$ never vanishes. Then we can consider the 
function on the boundary
\begin{equation}
\tau(x) = \frac{\partial a}{\partial x_3}(x)\cdot\frac{1}{a(x)}.
\end{equation}
Two natural questions now enter ones mind: is this $\tau$ globally 
well-defined, and what does this have to do with mean curvature?
\begin{prop}
Consider the immersion $\iota:\partial M\to M$.  Then the mean curvature $H$ 
satisfies
\begin{equation}
H = \frac{1}{2}\tau.
\end{equation}
\end{prop}
\noindent Since mean curvature is globally defined (on $\partial M$), so is $\tau$.
\begin{proof}
Let $\{x_1, x_2, x_3\}$ be Type A coordinates at the boundary, and let $g_{ij}$ and $\Gamma^m_{ij}$ be the corresponding metric tensor and Christoffel symbols, respectively.  Since the connection we are considering is the Levi-Civita connection, we have
\begin{equation*}
\Gamma^m_{ij}=\frac{1}{2}\sum_k \left(\frac{\partial}{\partial x_i}(g_{jk})+\frac{\partial}{\partial x_j}(g_{ik})- \frac{\partial}{\partial x_k}(g_{ij})\right)g^{km}.
\end{equation*}
(see, for example, \cite{doCarmo}).  By our choice of coordinate system, we have $g_{13}=g_{23}=g^{13}=g^{23}=0$ on the boundary and $\frac{\partial}{\partial x_1},\frac{\partial}{\partial x_2}$ are tangent to the boundary.  So on the boundary,
\begin{eqnarray}
\Gamma^1_{13}&=&\frac{1}{2}\sum_k \left(\frac{\partial}{\partial x_1}(g_{3k})+\frac{\partial}{\partial x_3}(g_{1k})- \frac{\partial}{\partial x_k}(g_{13})\right)g^{k1}\\
&=& \frac{1}{2} \left(\frac{\partial}{\partial x_1}(g_{31})+\frac{\partial}{\partial x_3}(g_{11})- \frac{\partial}{\partial x_1}(g_{13})\right)g^{11} + \\
& & \frac{1}{2} \left(\frac{\partial}{\partial x_1}(g_{32})+\frac{\partial}{\partial x_3}(g_{12})- \frac{\partial}{\partial x_2}(g_{13})\right)g^{21} \\
&=& \frac{1}{2}\left(\frac{\partial}{\partial x_3}(g_{11})g^{11}+\frac{\partial}{\partial x_3}(g_{12})g^{21}\right)\label{tau1}.
\end{eqnarray}
A similar calculation yields
\begin{eqnarray}\label{tau2}
\Gamma_{23}^2 = \frac{1}{2}\left(\frac{\partial}{\partial x_3}(g_{12})g^{12}+\frac{\partial}{\partial x_3}(g_{22})g^{22}\right).
\end{eqnarray}
Now, on the boundary, note that
\begin{eqnarray}\label{tau3}
g^{11}=\frac{g_{22}}{\det(g_{ij})},\mbox{ } g^{22}=\frac{g_{11}}{\det(g_{ij})},\mbox{ } g^{12}=g^{21}=-\frac{g_{12}}{\det(g_{ij})}.
\end{eqnarray}
Using a Laplace expansion on the bottom row of $(g_{ij})$, we also have
\begin{equation*}
\det(g_{ij}) = g_{31}\cdot\left|\begin{array}{cc} g_{12} & g_{13}\\ g_{22} &
g_{23}\end{array}\right| - g_{32}\cdot\left|\begin{array}{cc} g_{11} & g_{13}
\\ g_{21} & g_{23}\end{array}\right| + g_{33} \cdot\left|\begin{array}{cc} 
g_{11} & g_{12}\\ g_{21} & g_{22}\end{array}\right|.
\end{equation*}
On the boundary,
\begin{equation*}
g_{31}=\left|\begin{array}{cc} g_{12} & g_{13}\\ g_{22} &
g_{23}\end{array}\right| = g_{32} = \left|\begin{array}{cc} g_{11} & g_{13}
\\ g_{21} & g_{23}\end{array}\right| = 0.
\end{equation*}
Since $g_{33}\equiv 1$ everywhere, we have by product rule
\begin{equation}\label{tau4}
\frac{\partial}{\partial x_3}\det(g_{ij}) = \frac{\partial}{\partial x_3}\left|
\begin{array}{cc} g_{11} & g_{12}\\ g_{21} & g_{22}\end{array}\right| 
\mbox{  on $\partial M$.}
\end{equation}
Combining (\ref{tau1}), (\ref{tau2}), (\ref{tau3}), and (\ref{tau4}), we have
\begin{eqnarray*}
\Gamma_{13}^1 + \Gamma_{23}^2 &=& \frac{1}{2\det(g_{ij})}\left(\frac{\partial}{\partial x_3}(g_{11})g_{22}+\frac{\partial}{\partial x_3}(g_{22})g_{11}-2\frac{\partial}{\partial x_3}(g_{12})g_{12}\right)\\
&=&\frac{1}{2\det(g_{ij})}\frac{\partial}{\partial x_3}\left|\begin{array}{cc} g_{11} & g_{12}\\ g_{12} & g_{22}\end{array}\right|\\
&=& \frac{1}{2\det(g_{ij})}\frac{\partial}{\partial x_3}(\det(g_{ij}))\\
&=& \tau.
\end{eqnarray*}
Hence,
\begin{eqnarray*}
H &=& \frac{1}{2} \tr(S_{\nu}) = -\frac{1}{2}\left((\nabla_{\frac{\partial}{\partial x_1}}(-\frac{\partial}{\partial x_3}))_1 + (\nabla_{\frac{\partial}{\partial x_2}}(-\frac{\partial}{\partial x_3}))_2\right) \\
&=& \frac{1}{2}(\Gamma_{13}^1 + \Gamma_{23}^2) \\
&=& \frac{1}{2}\tau,
\end{eqnarray*}
as desired.
\end{proof}
We can also write $\tau$ in terms of Type B coordinates:
\begin{lem}\label{tauy}
Let $\{y_1, y_2, y_3\}$ be Type B coordinates.  Let $\{h_{ij}\}$ be 
the associated metric tensor and let $c:=\sqrt{\det(h_{ij})}$.  Then
\begin{equation*}
2\tau = 2{\sqrt{h_{33}}}\cdot{d(\sqrt{h_{33}}^{-1}(\nu))} + 2 \frac{dc(\nu)}{c}.
\end{equation*}
\end{lem}
\begin{proof}
Let $\{y_1,y_2,y_3\}$ be Type B coordinates.  Also consider Type A coordinates 
$\{x_1,x_2,x_3\}$ constructed with $u=y_3$.  Then 
$x_i\equiv y_i$ on the boundary for $i=1,2,3$.  In particular, this means that 
\begin{equation*}
b_{ij}:=\frac{\partial x_i}{\partial y_j} = \delta_{ij}\mbox{ for 
$i=1,2,3,j=1,2$},
\end{equation*}
where $B$ is the derivative matrix of the coordinate change.  Also, on the 
boundary we have 
\begin{equation*}
\frac{\partial}{\partial y_3} = (\sqrt{h_{33}})\nu = \sqrt{h_{33}}
\frac{\partial}{\partial x_3}.
\end{equation*}
So $b_{i3}=\sqrt{h_{33}}\delta_{i3}$. So $B$ is determined on $\partial M$. Also, on the 
boundary for $i,j=1,2$
\begin{eqnarray*}
0&=&\frac{\partial}{\partial y_i}(b_{j3})= \frac{\partial}{\partial y_i}
\left(\frac{\partial x_j}{\partial y_3}\right) \\
&=& \frac{\partial}{\partial y_3}\left(\frac{\partial x_j}{\partial y_i}
\right)=\frac{\partial}{\partial y_3}(b_{ji}).
\end{eqnarray*}
So on the boundary
\begin{equation*}
\frac{\partial}{\partial x_3}(b_{ij})\mbox{ for $i,j=1,2$}.
\end{equation*}
Hence, on the boundary
\begin{eqnarray*}
\det(B) = b_{33}=\sqrt{h_{33}},
\end{eqnarray*}
and
\begin{eqnarray*}
\frac{\partial}{\partial x_3}(\det(B))&=&\frac{\partial}{\partial x_3}
\left(b_{31}\cdot\left|\begin{array}{cc} b_{12} & b_{13}\\ b_{22} &
b_{23}\end{array}\right| - b_{32}\cdot\left|\begin{array}{cc} b_{11} & b_{13}
\\ b_{21} & b_{23}\end{array}\right| + b_{33} \cdot\left|\begin{array}{cc} 
b_{11} & b_{12}\\ b_{21} & b_{22}\end{array}\right|\right)\\
&=& \frac{\partial}{\partial x_3}(b_{33}) = \frac{\partial}{\partial x_3}
(\sqrt{h_{33}}).
\end{eqnarray*}
Let $g_{ij}$ be the metric tensor of the $x_i$'s and $a:=\sqrt{\det(g_{ij})}$.
Then
\begin{eqnarray*}
2\tau &=& 2 \frac{\partial a}{\partial x_3}\frac{1}{a} =
\frac{\partial}{\partial x_3}(\det(g_{ij}))\frac{1}{\det(g_{ij})} \\
&=& \frac{\partial}{\partial x_3}(\det(B^{-1})^2\det(h_{ij}))\frac{1}
{\det(B^{-1})^2\det(h_{ij})} \\
&=& (2\det(B^{-1})\frac{\partial}{\partial x_3}(\det(B^{-1}))\det(h_{ij})+\\
& & \det(B^{-1})^2\frac{\partial}{\partial x_3}(\det(h_{ij})))
\frac{1}{\det(B^{-1})^2\det(h_{ij})} \\
&=& 2\sqrt{h_{33}}\frac{\partial}{\partial x_3}(\sqrt{h_{33}}^{-1}) +
\frac{\partial}{\partial x_3}(\det(h_{ij}))\frac{1}{\det(h_{ij})}.
\end{eqnarray*}
Since
\begin{equation*}
\frac{\partial}{\partial x_3}(\det(h_{ij}))\frac{1}{\det(h_{ij})} = 
2\frac{\partial c}{\partial x_3}\frac{1}{c}
\end{equation*}
and $\frac{\partial}{\partial x_3}=\nu$, we have the result.
\end{proof}
\end{section}
%
%
%
%
\begin{section}{The Image of the Curvature Form}
We will use this $\tau$ to prove the following lemma, which relates to the 
image of the curvature form.
\begin{lem}\label{bct}
Suppose $M$ is a $3$-manifold with boundary, $k+1>3/2$, $\alpha,\beta\in H^{k+1}_{con}(\kb)\cap H_A$ and $\nabla_{A}\in\Ck$.  Then 
\begin{equation}\label{thecondition1}
d_A[\alpha\cdot\beta](\nu)= -2\tau[\alpha\cdot\beta]\mbox{    on $\partial M$,}
\end{equation}
where $\nu$ is the normal inward pointing vector field.
\end{lem}

\noindent Since $k+1>3/2$, note that $[\alpha\cdot\beta]$ is $C^1$, and thus $d_A[\alpha\cdot\beta]$ 
is continuous.  Hence, the above equality is true not just in the trace sense, but as 
an equality of two continuous functions.
\begin{proof}
We will use Type A coordinates $(x_1,x_2,x_3)$, and assume 
that the vector bundle $\k_P$ is also trivialized in this neighborhood.  
Recall that the metric tensor in this coordinate system has the feature that 
$g_{i3}=\delta_{i3}$ on the boundary, and $g_{33}=1$ everywhere.  Thus, 
$g^{i3}=\delta_{i3}$ also on the boundary.  Also, $\frac{\partial}
{\partial x_3}$ is the inward pointing 
normal vector on the boundary.  
Take $\alpha,\beta$ as above and define 
$\alpha_i$ and $\beta_i$ so that $\alpha = \sum_{i=1}^3\alpha_i dx_i$ and 
$\beta = \sum_{i=1}^3\beta_i dx_i$.  Since we are assuming $\k_P$ has a fixed 
trivialization in our neighborhood, we can view the $\alpha_i$ and $\beta_i$ 
as $\k$-valued functions.
Also, since $\frac{\partial}{\partial x_3}$ 
is the inward pointing normal vector and $\alpha,\beta$ satisfy conductor 
boundary conditions, we have 
\begin{equation}\label{bc1}
\alpha_1=\alpha_2=\beta_1=\beta_2=0\mbox{   on $\partial M$.}
\end{equation}
Let $d$ be the flat connection with respect to our fixed trivialization of 
$\k_P$ and define a $\k$-valued $1$-form $A$ 
so that $d_A = d + A$.  Define $A_i$ so that $A=\sum_{i=1}^3 A_idx_i$. On 
this coordinate patch, we have
\begin{eqnarray}\label{bc2}
[\alpha\cdot\beta] &=& \sum_{j,k=1}^3 [\alpha_j,\beta_k](dx_j\cdot dx_k)
= \sum_{j,k=1}^3 [\alpha_j,\beta_k]g^{jk}.
\end{eqnarray}
Taking the derivative $d_A$ yields
\begin{eqnarray*}
d_A([\alpha\cdot\beta])&=& \sum_{j,k=1}^3 d_A([\alpha_j,\beta_k]g^{jk})\\
&=& \sum_{j,k=1}^3 [d_A(\alpha_j),\beta_k]g^{jk} + [\alpha_j,d_A(\beta_k)]g^{jk} + [\alpha_j,\beta_k]d(g^{jk})
\end{eqnarray*}
If both $j,k<3$, then $\alpha_j=\beta_k=0$, and hence 
\begin{equation}\label{bc1.5}
[d_A(\alpha_j),\beta_k]g^{jk} + [\alpha_j,d_A(\beta_k)]g^{jk} + [\alpha_j,\beta_k]d(g^{jk}) = 0\mbox{ on $\partial M$}.
\end{equation}
Suppose $j=3$ and $k<3$.  Then $\beta_k=0$ on $\partial M$ by \eqref{bc1}, and $g^{3k}=0$ on $\partial M$ since we are using 
Type A coordinates.  Thus, \eqref{bc1.5} holds in this case also.
Similiarly, if $j< 3$ and $k=3$ then $\alpha_j=0$ and $g^{j3}=0$ and thus \eqref{bc1.5} holds.  
In sum, we have
\begin{eqnarray}\label{bc2.5}
d_A([\alpha\cdot\beta])|_{\partial M} = [d_A\alpha_3,\beta_3]+
[\alpha_3,d_A\beta_3] + [\alpha_3,\beta_3]d(g^{33}).
\end{eqnarray}
Using the adjoint matrix, we see that
\begin{equation*}
g^{33}=\frac{\left|\begin{array}{cc} g_{11} & g_{12}\\ g_{21} & g_{22}\end{array}\right|}{\det(g_{ij})}.
\end{equation*}
Combining the fact that $\left|\begin{array}{cc} g_{11} & g_{12}\\ g_{21} & g_{22}\end{array}\right|=\det(g_{ij})$ on 
$\partial M$ and \eqref{tau4} yields
\begin{eqnarray}\label{bc2.75}
\frac{\partial g^{33}}{\partial x_3} &=& \frac{\left(\frac{\partial}{\partial x_3} \left|
\begin{array}{cc} g_{11} & g_{12}\\ g_{21} & g_{22}\end{array}\right|\right)\det(g_{ij})-
\left|
\begin{array}{cc} g_{11} & g_{12}\\ g_{21} & g_{22}\end{array}\right|(\frac{\partial}{\partial x_3}\det(g_{ij}))}{\det(g_{ij})^2}\\
&=& 0 \mbox{ on $\partial M$.}
\end{eqnarray}
Hence, we have
\begin{equation}\label{bc3}
d_A([\alpha\cdot\beta])|_{\partial M}(\frac{\partial}{\partial x_3}) = [d_A\alpha_3(\frac{\partial}{\partial x_3}),\beta_3]+
[\alpha_3,d_A\beta_3(\frac{\partial}{\partial x_3})].
\end{equation}

We will leave $d_A([\alpha\cdot\beta])|_{\partial M}$ for the moment 
and investigate what $d^{*}_A\alpha = d^{*}_A\beta = 0$ 
means in our coordinate system.  We will calculate $d^{*}$ by using the 
Hodge star operator.  One can verify that
\begin{eqnarray*}
*dx_j &=& a(g^{j1}dx_2\wedge dx_3 + g^{j2}dx_3\wedge dx_1 + g^{j3}dx_1
\wedge dx_2), \\
\end{eqnarray*}
where $a=\sqrt{\det(g_{ij})}=1/\sqrt{\det(g^{ij})}$.  Using this, we calculate
\begin{eqnarray*}
-d^* \alpha &=& *d*(\alpha_1 dx_1 + \alpha_2 dx_2 + \alpha_3 dx_3) \\
&=& *d\left(a\left(\sum_j \alpha_j(g^{j1}dx_2\wedge dx_3 + g^{j2}dx_3\wedge dx_1 + g^{j3}dx_1
\wedge dx_2)\right)\right) \\
&=& *(da)\left(\sum_j \alpha_jg^{j1}dx_2\wedge dx_3 + \alpha_jg^{j2}dx_3\wedge dx_1 + \alpha_jg^{j3}dx_1
\wedge dx_2\right) + \\
& & \left(\sum_j \frac{\partial}{\partial x_1}(\alpha_jg^{j1}) + \frac{\partial}{\partial x_2}(\alpha_jg^{j2}) + \frac{\partial}{\partial x_3}(\alpha_jg^{j3})\right)*(adx_1\wedge dx_2\wedge dx_3) \\
&=& *(da)\left(\sum_j \alpha_jg^{j1}dx_2\wedge dx_3 + \alpha_jg^{j2}dx_3\wedge dx_1 + \alpha_jg^{j3}dx_1
\wedge dx_2\right) + \\
& & \left(\sum_{i,j} \frac{\partial}{\partial x_i}(\alpha_jg^{ji})\right).
\end{eqnarray*}
By (\ref{bc1}) and the fact that $\frac{\partial}{\partial x_1},
\frac{\partial}{\partial x_2}$ are tangent at the boundary, we have that 
for $i=1,2,j=1,2$
\begin{eqnarray*}
\frac{\partial}{\partial x_i}(\alpha_jg^{ji})|_{\partial M} &=& 
(\frac{\partial}{\partial x_i}(\alpha_j)g^{ji} + \alpha_j\frac{\partial}
{\partial x_i}(g^{ji}))|_{\partial M}  \\
&=& 0\cdot g^{ji} + 0\cdot\frac{\partial}{\partial x_i}(g^{ji}) = 0.
\end{eqnarray*}
Also, since $\frac{\partial}{\partial x_1},\frac{\partial}{\partial x_2}$ are 
tangent to the boundary, we also have for $i=1,2$
\begin{eqnarray*}
\frac{\partial}{\partial x_i}(\alpha_3g^{3i})|_{\partial M} &=& 
(\frac{\partial}{\partial x_i}(\alpha_3)g^{3i} + \alpha_3\frac{\partial}
{\partial x_i}(g^{3i}))|_{\partial M}  \\
&=& \frac{\partial}{\partial x_i}(\alpha_3)\cdot 0 + \alpha_3\cdot 0 = 0.
\end{eqnarray*}
Finally, for $j=1,2$:
\begin{eqnarray*}
\frac{\partial}{\partial x_3}(\alpha_jg^{j3})|_{\partial M} &=& 
(\frac{\partial}{\partial x_3}(\alpha_j)g^{j3} + \alpha_j\frac{\partial}
{\partial x_3}(g^{j3}))|_{\partial M}  \\
&=& \frac{\partial}{\partial x_3}(\alpha_j)\cdot 0 + 0\cdot\frac{\partial}
{\partial x_3}(g^{j3}) = 0.
\end{eqnarray*}
Hence, $-d^*\alpha$ on the boundary reduces to
\begin{eqnarray*}
-d^*\alpha|_{\partial M} &=&  \frac{\partial}{\partial x_3}(g^{33}\alpha_3) + 
*(da)(g^{33}\alpha_3 dx_1\wedge dx_2).
\end{eqnarray*}
We showed that $\partial/\partial x_3(g^{33})=0$ in \eqref{bc2.75}.
Thus, we obtain
\begin{eqnarray*}
-d^*\alpha|_{\partial M} &=&  \frac{\partial}{\partial x_3}(\alpha_3) + 
*(da)(g^{33}\alpha_3 dx_1\wedge dx_2) \\
&=& \frac{\partial}{\partial x_3}(\alpha_3) + \frac{1}{a}\frac{\partial a}
{\partial x_3}(\alpha_3).
\end{eqnarray*}
Since $d_A^*\alpha =d^*\alpha - [A\cdot\alpha]$, we have,
\begin{eqnarray}\label{bc4}
-d^*_A\alpha|_{\partial M}= \frac{\partial \alpha_3}{\partial x_3} + 
\frac{1}{a}\frac{\partial a}{\partial x_3}(\alpha_3) - [\alpha\cdot A]\\
=\frac{\partial \alpha_3}{\partial x_3} + 
\frac{1}{a}\frac{\partial a}{\partial x_3}(\alpha_3) + [A_3,\alpha_3]
\label{bc5},
\end{eqnarray}
where we used (\ref{bc1}) and (\ref{bc2}) (replacing $\beta$ with $A$) in the 
last line.  Of course, an analogous statement holds for $\beta$ replacing $\alpha$.  

We now revisit (\ref{bc3}) and plug in (\ref{bc5}):
\begin{eqnarray*}
d_A([\alpha\cdot\beta])|_{\partial M}(\nu) &=& 
[d_A\alpha_3,\beta_3](\partial/\partial x_3)+
[\alpha_3,d_A\beta_3](\partial/\partial x_3) \\
&=& [\frac{\partial\alpha_3}{\partial x_3},\beta_3]+[[A_3,\alpha_3],\beta_3]
\\ & & + [\alpha_3,\frac{\partial\beta_3}{\partial x_3}] + 
[\alpha_3,[A_3,\beta_3]] \\
&=& -[\frac{1}{a}\frac{\partial a}{\partial x_3}\alpha_3 + [A_3,\alpha_3]+d_A^*\alpha,
\beta_3]+[[A_3,\alpha_3],\beta_3] \\
& & - [\alpha_3,\frac{1}{a}\frac{\partial a}
{\partial x_3}\beta_3 + [A_3,\beta_3]+d_A^*\beta] + [\alpha_3,[A_3,\beta_3]] \\
&=& -\frac{2}{a}\frac{\partial a}{\partial x_3}[\alpha_3,\beta_3]-[d^*_A\alpha,\beta_3]-[\alpha_3, d^*_A\beta] \\
&=& -2\tau[\alpha\cdot\beta]|_{\partial M}-[d^*_A\alpha,\beta(\nu)]-[\alpha(\nu), d^*_A\beta],\\
&=& -2\tau[\alpha\cdot\beta]|_{\partial M}.
\end{eqnarray*}
where we again used (\ref{bc2}) on the second to last line, as well as the fact that $\alpha,\beta\in H_A$.  The lemma is thus proven.
\end{proof}
Inspired by the previous result, we define a linear map 
\begin{equation*}
T_A: \Lie(\Gk)\to L^2(\kb|_{\partial M})
\end{equation*}
given by
\begin{equation}
T_A(f)= d_A\Delta_A f + 2\tau\Delta_A f.
\end{equation}
Counting derivatives (note that $\Lie(\Gk)=H^{k+1}_{con}(\kb)$), and using Theorem 9.3 in 
\cite{Palais}, we see that $T_A$ is well-defined and bounded. Define a set $\mathcal{L}_A\subseteq\Lie(\Gk)$ as 
\begin{equation}
\mathcal{L}_A:=\Span\{\mathcal{R}_A(\alpha\cdot\beta):\alpha,\beta\in H_A\}
\end{equation}
The previous lemma yields
\begin{cor}
The set $\mathcal{L}_A$ is contained in $\mathrm{ker}(T_A)$.  In particular, since 
$T_A$ is not identically $0$, we have that $\overline{\mathcal{L_A}}$ is a proper 
subset of $\Lie(\Gk)$.
\end{cor}
\begin{proof}
Let $g\in\mathcal{L}_A$.  By Lemma \ref{bct}, since 
\begin{equation*}
\mathcal{R}_A(\alpha,\beta) = -2G_A[\alpha\cdot\beta],
\end{equation*}
$\Delta_A f$ satisfies
\begin{equation*}
T_A(g)=d_A(\Delta_A g)(\nu) + 2\tau (\Delta_A g) = 0,
\end{equation*}
proving the corollary.
\end{proof}
This corollary shows that the image of the curvature form can never be dense in the 
gauge algebra, unlike the case in \cite{NR}.
\end{section}
%
%
%
%
\begin{section}{A Partial Converse of Lemma \ref{bct}}
The natural question now is whether the converse of Lemma \ref{bct} holds.  A 
quick argument shows that it cannot.  
\begin{lem}
The converse to Lemma \ref{bct} does not hold.  More specifically, 
there exists $f\in\mathrm{ker}(T_A)-\mathcal{L_A}$ if the connection $\nabla_{A}\in\Ck$ 
also lies in $\mathcal{C}^{k+1}_{con, A_0}$.
\end{lem}
\begin{proof}
Suppose $\nabla_A\in\mathcal{C}^{k+1}_{con, A_0}$.
Note that $\mathcal{L}_A\subseteq H^{k+2}_{con}(\kb)$ since $\nabla_{A}\in\mathcal{C}^{k+1}_{con, A_0}$, 
and thus $G_A: H^k(\kb)\to H^{k+2}_{con}(\kb)$ exists by Proposition \ref{Greenexist}.  However, the domain of $T_A$ is 
$\Lie(\Gk)=H^{k+1}_{con}(\kb)$.  This disparity in regularity will sink the converse as 
follows.  We first construct $f\in H^{k-1}(\kb)$ that is not in
$H^k(\kb)$ and is $0$ in a neighborhood of the boundary:  Take 
an open subsets $U\subset W\subset M$ such that $\bar{U}\subset W$, $\bar{W}\subset \mathring{M}$
and take $\tilde{f}\in H^{k-1}(\kb|_U) - H^{k}(\kb|_U)$ and $\tilde{f}\in H^{k-1}(\bar{W})$. Take a smooth function 
$\zeta:M\to [0,1]$ such that $\zeta|_{\bar{U}} \equiv 1$ and $supp(\zeta)\subset W$.  Then 
$f:=\zeta\cdot \tilde{f}\in H^{k-1}(\kb|_{W})$ , and we can extend $f$ by $0$ to have 
$f\in H^{k-1}(\kb)$.  By the equivalence in Section 4 of \cite{Palais}, 
since $f|_U=\tilde{f}\notin H^{k}(\kb|_U)$ we have $f\notin H^{k}(\kb)$.  However, we do 
have that $f$ is 0 in a neighborhood of the boundary.  Set $g:=G_A f$.  Then 
$g\in H^{k+1}_{con}(\kb)=\Lie(\Gk)$.  If $g\in H^{k+2}(\kb)$, then $f=\Delta_A g\in H^k$ which 
is a contradiction.  So $g\in H^{k+1}_{con}(\kb)-H^{k+2}_{con}(\kb)$, and $\Delta_A g = f$ vanishes in a neighborhood of 
$\partial M$, so $d_A(\Delta_A f)(\nu) = 0 = 2\tau(\Delta_A f)$.  So, $g\in\mathrm{ker}(T_A)$, but not in $H^{k+2}(\kb)$ 
and thus not in $\mathcal{L}_A$, as desired.
\end{proof}

In the previous lemma, regularity considerations sunk the converse.  However, if we took regularity 
out of the equation, perhaps the converse would hold.  In other words, perhaps we have
\begin{equation*}
\mathrm{ker}T_A\cap C^\infty = \mathcal{L}_A\cap C^\infty.
\end{equation*}
So that the Green operator $G_A$ maps smooth functions to smooth functions, we also want 
the connection $\nabla_{A}$ to be $C^\infty$.  In this setting, the converse does for a 
specific set up.  Namely, if $P$ is the trivial bundle 
$\bar{O}\times K\to \bar{O}$ for a bounded open set $O\subseteq\mathbb{R}^3$ with smooth boundary with and the base connection $\nabla_{A_0}$ is the flat connection.  In this set up, 
$\Kb$ is isomorphic to $\bar{O}\times K\to \bar{O}$, and $\kb$ is isomorphic to 
$\bar{O}\times\mathfrak{k}\to \bar{O}$.  So we can view gauge transformations $g$ as $K$-valued 
functions on $\bar{O}$, gauge algebra elements $\psi$ as $\mathfrak{k}$-valued 
functions, and $\kb$-valued forms as $\mathfrak{k}$-valued forms.

As in Chapter 2, we denote the flat connection as $\nabla_0$.  This means we should denote exterior differentiation
by $d_0$, but since $\nabla_0 = d$ (as asserted in Chapter 2), we will instead simply use $d$ without a subscript.  Similiarly, 
we denote $d^*_0$ by simply $d^*$. 

To do prove the converse of Lemma \ref{bct} in this case, we proceed locally.
We first consider interior neighborhoods.  For this, we prove a lemma originally 
outlined by L. Gross.

\begin{lem}\label{noboundary1}
Let $(a,b)^3$ be a cube in $\mathbb{R}^3$, and let $\Psi:(a,b)^3\to\k$ be 
smooth and 
have compact support.  Then $\Psi\in \Span\{[\alpha\cdot\beta]: \alpha,\beta
\in C^\infty_c(\Lambda^1((a,b)^3\otimes\k)), d^*\alpha =d^*\beta = 0\}$.
\end{lem}
\begin{proof}
Since $\k=[\k,\k]$ by semisimplicity, without loss of generality we can assume 
$\Psi(x,y,z) = \psi(x,y,z)\cdot[A,B]$, where 
$A,B$ are fixed elements of $\k$ and $\psi\in C^\infty_c((a,b)^3)$.  Choose $c,
d,i,j,k,l\in\mathbb{R}$ so that $supp(\psi)\subset (c,d)^3$ and 
$a<j<k<i<c<d<l<b$.  We can 
then find a function $h:(a,b)^3\to\mathbb{R}$ so that 
\begin{enumerate}
\item\label{p1} $h$ is smooth
\item\label{p2}$h|_{[c,d]^3}\equiv \psi$, 
\item\label{p3}$supp(h)\subset(j,d)^3$,
\item\label{p4}$\int_j^d h(x,s,z)ds = 0$ for any fixed $x,z$.
\item\label{p5} $h(x,y,z) = 0$ if $(x,z) \notin (c,d)^2$.
\item\label{p6} $h|_{y\in (i,c)}\equiv 0$.
\end{enumerate}
Specifically, define 
\begin{eqnarray*}
\eta(t)= \left\{\begin{array}{l@{\quad\mbox{ if }\quad}l}C\exp\left(\frac{1}
{1-(\frac{2}{i-k}(t - 
\frac{i+k}{2}))^2}\right) & t\in(k,i) \\
0 & \mbox{  else}\end{array}\right.
\end{eqnarray*}
where $C$ is chosen so the integral of $\eta$ is $1$. Let $I(x,z):=\int_c^d
\psi(x,s,z)ds$, and finally define $h(x,y,z):=-I(x,z)\eta(y)+\psi(x,y,z)$.  
One can check that $h$ satisfies the above properties.

Define $F:(a,b)^3\to \mathbb{R}$ as $F(x,y,z)=\int_a^y h(x,s,z)ds$.  Then $F$ 
is smooth and 
$supp(f)\subset(j,d)^3$ 
by Properties \ref{p3} and \ref{p4} above.  Also, it is clear that $F_y = h$.  
We now construct another function $G:(a,b)^3\to\mathbb{R}$.  Define $G$ as 
$G(x,y,z)=\phi(x)v(y,z)$, where $\phi:(a,b)\to\mathbb{R}$ is constructed so 
$\phi\in C^\infty_c(a,b)$ and $\phi|_{[c,d]}(x)=x$, and $v:(a,b)^2\to
\mathbb{R}$ is constructed so $v|_{[c,d]^2}\equiv 1$, and 
$supp(v)\subset (i,l)^2$.  Then $G_x|_{[c,d]^3}\equiv 1$, and has compact 
support.  Let us consider 
$\Theta(x,y,z):=F_y(x,y,z)\cdot G_x(x,y,z) = h(x,y,z)\cdot G_x(x,y,z)$. We 
will show that $\Theta(x,y,z)=\psi(x,y,z)$ by looking at it in cases:  

First, if $(x,z)\notin (c,d)^2$, then by property \ref{p5} we have 
$h(x,y,z)=0=\psi(x,y,z)$.  So for our next cases we can assume $(x,z)\in(c,d)^2$, and 
thus $G_x(x,y,z)=v(y,z)$.  If $y\in(a,i]$, then $v(y,z)=0$ since 
$supp(v)\subset (i,l)^2$.  Hence, 
$\Theta(x,y,z)=0=\psi(x,y,z)$.  If $y\in (i,c)$, then $h(x,y,z)=0$ by 
property \ref{p6}, and so $\Theta(x,y,z)=0=\psi(x,y,z)$.  
If $(x,y,z)\in [c,d]^3$, then $\Theta(x,y,z) = \psi(x,y,z)\cdot 1 = 
\psi(x,y,z)$.  And finally, if $y\in (d,b)$, then $h(x,y,z)=0=\psi(x,y,z)$.  
Hence, in all cases, $\Theta(x,y,z)=\psi(x,y,z)$.

Now define 2-forms $\omega_1, \omega_2$ as $\omega_1 = -F\cdot A dy\wedge dz$ 
and $\omega_2 = G\cdot B dz\wedge dx$.  Let $\alpha:=d^*\omega_1$ and 
$\beta:=d^*\omega_2$.  Since $(d^*)^2=0$, we have $d^*\alpha = d^*\beta =0$.  
Now,
\begin{eqnarray*}
\alpha &=& d^*(\omega_1) = *d*(\omega_1)=*d(-F\cdot Adx)\\
&=&*(F_y\cdot Adx\wedge dy - F_z \cdot A dz\wedge dx)B= - F_z\cdot A dy + 
F_y \cdot A dz.
\end{eqnarray*}
Similarly,
\begin{eqnarray*}
\beta &=& *d(G\cdot B dy)= *(G_x B dx\wedge dy - G_z B dy\wedge dz) = - G_z
\cdot B dx + G_x\cdot B dz
\end{eqnarray*}
Since $\Theta=\psi$,
\begin{eqnarray*}
[\alpha\cdot \beta] &=& [0,G_z\cdot B] + [F_z\cdot A,0] + [F_y A, G_x B] \\
&=& \psi[A,B]=\Psi,
\end{eqnarray*}
as desired.
\end{proof}
We extend this result to any domain $O$.

\begin{lem}\label{noboundary2}
Let $O\subset\mathbb{R}^3$ be a bounded open set.  Then $C^\infty_c(O\otimes\k)
=\Span\{[\alpha\cdot\beta]: \alpha,\beta\in C^\infty_c(\Lambda^1(O\otimes\k)), 
d^*\alpha =d^*\beta = 0\}$.
\end{lem}
\begin{proof}
Let $\Psi\in C^\infty_c(O\otimes\k)$ be arbitrary.  Let $\{C_k\}_{k=1}^n$ be a 
finite family of open cubes that cover the support of $\Psi$ and are contained
in $O$.  Let $\{\lambda_k\}$ be a partition of unity subordinate to 
the cover $\{C_k\}$.  Then the function $\lambda_k\cdot\Psi$ lies in $C^\infty
_c(C_k\otimes\k)$.  By the previous lemma, there exists sequences 
$\{\alpha_i\}_{i=1}^{i(k)}, \{\beta_i\}_{i=1}^{i(k)}$ such that each 
$\alpha_i,\beta_i\in C^\infty_c(\Lambda^1(C_k\otimes\k))$, $d^*\alpha_i=
d^*\beta_i=0$ and 
$\lambda_k\cdot\Phi = \sum_{i=1}^{i(k)}[\alpha_i\cdot\beta_i]$ on $C_k$.  
Extending the $\alpha$'s and $\beta$'s by zero, we have 
$\alpha_i,\beta_i\in C^\infty_c
(\Lambda^1(O\otimes\k))$, $d^*\alpha=d^*\beta=0$, and $\lambda_k\cdot\Phi = 
\sum_{i=1}^{i(k)}[\alpha_i\cdot\beta_i]$ on $0$.  Thus, $\lambda_k\cdot\Phi\in
Span\{[\alpha\cdot\beta]: \alpha,\beta\in C^\infty_c(\Lambda^1(O\otimes\k)), 
d^*\alpha =d^*\beta = 0\}$.  So, $\Psi=\sum_{k=1}^n (\lambda_k\cdot\Phi)\in
Span\{[\alpha\cdot\beta]: \alpha,\beta\in C^\infty_c(\Lambda^1(O\otimes\k)), 
d^*\alpha =d^*\beta = 0\}$, as desired.
\end{proof}

%
%
We now look at neighborhoods of the boundary of $O$ and see if all the smooth $\Psi$ that 
satisfy the boundary condition of Lemma \ref{bct} are in the desired span.  We 
see that this is so.
\begin{lem}\label{lbo}
Let $O\subset\mathbb{R}^3$ be open and bounded, and let $U$ be a neighborhood 
of $\bar{O}$ that includes the boundary, admits the Type B coordinates  
$\{y_1,y_2,y_3\}$, and
is a cube under these coordinates.
Let $\Psi:U\to\k$ be smooth, have compact support, and  
\begin{equation}\label{bco1}
d\Psi(\nu) = -2 \tau \Psi\mbox{  on $\partial O\cap U$.}
\end{equation}
Then $\Psi\in \Span\{[\alpha\cdot\beta]: \alpha,\beta
\in C^\infty_c(U\otimes\k)); d^*\alpha =d^*\beta = 0; \mbox{$\alpha,\beta$
satisfy CBC}\}$.
\end{lem}
In the preceeding lemma and in what follows, a smooth $1$-form $\alpha$ satisfies 
\emph{conductor 
boundary conditions} (or \emph{CBC} for short) if $\alpha\in H^1_{con}(\kb)$ as well as 
being smooth.  This is equivalent to saying that $\iota^*(\alpha)=0$, where 
$\iota:\partial M\to M$ is the inclusion, or saying the tangential component of $\alpha$ on
the boundary is $0$.  Also, viewing $U$ as the cube $(0,1)\times (0,1)
\times [0,1)$, a function $f\in C^\infty_c(U)$ has its support contained in $(\epsilon,
1-\epsilon)\times (\epsilon,1-\epsilon)\times [0,1-\epsilon)$ for some 
$\epsilon>0$.  The point is that it need not vanish on the boundary 
$\{y_3=0\}$.
\begin{proof}[Proof of Lemma \ref{lbo}]
Let $\{v_i\}$ be basis of $\k$.  Then we can write 
$\Psi = \sum \psi_i\cdot v_i$.  Since the basis elements are independent, 
by (\ref{bco1}) we have that $d\psi_i(\nu)=-2\tau\psi_i$.  Since 
$\mathfrak{k}$ is semisimple, each basis element $v_i$ can be written as a sum 
of commutators $v_i=\sum_{j=1}^{\alpha(i)}[f_j^i,g_j^i]$.  Hence, we can write 
$\Psi$ as
\begin{equation*}
\Psi = \sum_i \sum_{j=1}^{\alpha(i)}\psi_i [f_j^i,g_j^i].
\end{equation*}
So without loss of generality we can assume $\Psi = \psi\cdot[A,B]$, where 
$A,B$ are fixed elements of $\k$ and $\psi\in C^\infty_c(U)$ and 
$d\psi(\nu)=-2\tau\psi$.

Coordinatize $U$ using Type B coordinates $(y_1,y_2,y_3)$ under which 
the domain is a cube.  Then, without loss of generality, $U =  
(0,\delta)\times(0,\delta)\times[0,\delta)$ under these coordinates.  Here, let 
$a:=\sqrt{\det(h_{ij})}$, where $h_{ij}$ is the metric tensor of our chart.

Choose $c,d,i,j,k,l\in\mathbb{R}$ so that $supp(\psi)\subset(c,d)\times(c,d)
\times[0,d)$ and 
$0<j<k<i<c<d<l<\delta$.  We can define a function $h$ similar to the $h$ in Lemma 
\ref{noboundary1}. This time we set $h(y_1,y_2,y_3) = -I(y_1,y_3)\eta(y_2) + 
h^{33}(y_1,y_2,y_3)a(y_1,y_2,y_3)^2\psi(y_1,y_2,y_3)$, where $\eta$ is the 
same bump function 
from Lemma \ref{noboundary1}, and 
\begin{equation*}
I(y_1,y_3)=\int_c^d h^{33}(y_1,s,y_3)a(y_1,s,y_3)^2\psi(y_1,s,y_3)ds.
\end{equation*}  
Then $h$ has the following properties, 
analogous to the previous case in Lemma \ref{noboundary1}:
\begin{enumerate}
\item\label{p11} $h$ is smooth
\item\label{p12}$h|_{[c,d]\times[c,d]\times[0,d]}\equiv 
h^{33}(y_1,y_2,y_3)a(y_1,y_2,y_3)^2\psi|_{[c,d]\times[c,d]\times[0,d]}$, 
\item\label{p13}$supp(h)\subset(c,d)\times(j,d)\times[0,d)$,
\item\label{p14}$\int_0^{y_2} h(y_1,s,y_3)ds = 0$ for any fixed $y_1,y_3$ and 
$y_2\ge d$ or $y_2\le j$,
\item\label{p16} $h|_{(i,c)}\equiv 0$.
\end{enumerate}

$h$ has an additional property.  Using Lemma \ref{tauy}, we have
\begin{eqnarray}\label{bco3}
\frac{\partial (h^{33}a^2\cdot\psi)}{\partial y_3} &=& 
2\sqrt{h^{33}}\frac{\partial\sqrt{h^{33}}}{\partial y_3}a^2\psi +
h^{33} 2a\cdot\frac{\partial a}{\partial y_3}\cdot\psi +
h^{33}a^2\cdot\frac{\partial\psi}{\partial y_3}\\
&=& h^{33}a^2\left[2\left(\frac{1}{\sqrt{h^{33}}}\frac{\partial\sqrt{h^{33}}}
{\partial y_3}a^2 + \frac{1}{a}\cdot\frac{\partial a}{\partial y_3}\cdot
\right)\psi + \frac{\partial\psi}{\partial y_3}\right] \\
&=& 0 \mbox{  on $\partial O\cap U$}.
\end{eqnarray}
Hence, differentiating under the integral yields
\begin{equation}\label{bco4}
\frac{\partial h}{\partial y_3} = 0\mbox{  on $\partial O\cap U$}.
\end{equation}

Define $F:U\to \mathbb{R}$ as $F(y_1,y_2,y_3)=\int_0^{y_2} h(y_1,s,y_3)ds$.  
Then $F$ is smooth and $supp(F)\subset(c,d)\times(j,d)\times[0,d)$ 
be Properties \ref{p13} and \ref{p14} above.  Also, it is clear that 
$F_2 = h$.  Also, by (\ref{bco4}), differientiating under the integral sign 
yields
\begin{equation}\label{bco5}
\frac{\partial F}{\partial y_3} = 0\mbox{  on $\partial O\cap U$}.
\end{equation}

We now construct another function $G:[0,\delta]^3\to\mathbb{R}$ which is completely
analogous to the $G$ in Lemma \ref{noboundary1}.  Define $G$ as 
$G(y_1,y_2,y_3)=v_1(y_1)v_2(y_2)v_3(y_3)$, where $v_i:[0,1]\to\mathbb{R}$ is 
constructed as follows:
$v_1\in C^\infty_c(0,\delta)$, $v_1|_{[c,d]}(x)=x$,and $supp(v_1)\subset (i,l)$; 
$v_2\in C^\infty_c(0,\delta)$, $v_2|_{[c,d]}\equiv 1$, and $supp(v_2)\subset (i,l)$;
$v_3\in C^\infty_c([0,\delta))$, $v_3|_{[0,d]}\equiv 1$, and $supp(v_3)\subset 
[0,l)$.Then $G_1|_{[c,d]\times[c,d]\times[0,d]}
\equiv 1$, and has compact support in $U$.  Let us again consider 
$\Theta(y_1,y_2,y_3):=F_2(y_1,y_2,y_3)\cdot G_1(y_1,y_2,y_3) = h(y_1,y_2,y_3)
\cdot G_1(y_1,y_2,y_3)$, and show that $\Theta=h^{33}a^2\psi$ by looking at it 
in cases:  

First, if $(y_1,y_3)\notin (c,d)\times[0,d)$, then by property \ref{p13} we 
have $h(y_1,y_2,y_3)=0=h^{33}(y_1,y_2,y_3)a(y_1,y_2,y_3)^2\psi(y_1,y_2,y_3)$.  
So now we 
can assume $(y_1,y_3)\in(c,d)\times[0,d)$, and 
thus $G_1(y_1,y_2,y_3)=v_2(y_2)$.  If $y_2\in(0,i]$, then $v_2(y_2)=0$ 
since $supp(v_2)\subset (i,l)$.  Hence, 
$\Theta(y_1,y_2,y_3)=0=h^{33}(y_1,y_2,y_3)a(y_1,y_2,y_3)^2\psi(y_1,y_2,y_3)$.  
If $y_2\in (i,c)$, 
then $h(y_1,y_2,y_3)=0$ by property \ref{p16}, and so 
$\Theta(y_1,y_2,y_3)=0=h^{33}(y_1,y_2,y_3)a(y_1,y_2,y_3)^2\psi(y_1,y_2,y_3)$.  
If $y_2\in [c,d]$, then $(y_1,y_2,y_3)\in [c,d]^2\times[0,d]$. So 
\begin{eqnarray*}
\Theta(y_1,y_2,y_3) &=& h^{33}(y_1,y_2,y_3)a(y_1,y_2,y_3)^2\psi(y_1,y_2,y_3)
\cdot 1 \\
&=& h^{33}(y_1,y_2,y_3)a(y_1,y_2,y_3)^2\psi(y_1,y_2,y_3).
\end{eqnarray*}
And finally, if $y_2\in (d,\delta)$, then by 
property \ref{p13}
\begin{equation*} 
h(y_1,y_2,y_3)=0=h^{33}(y_1,y_2,y_3)a^2(y_1,y_2,y_3)\psi(y_1,y_2,y_3).
\end{equation*}  
Hence, in all cases, $\Theta(y_1,y_2,y_3)=h^{33}(y_1,y_2,y_3)a^2(y_1,y_2,y_3)
\psi(y_1,y_2,y_3)$.

Now define 2-forms $\omega_1, \omega_2$ as $\omega_1 =-F\cdot A(*^{-1}(dy_1))$ 
and $\omega_2 = G\cdot B (*^{-1}(dy_2))$.  Let $\alpha:=d^*\omega_1$ and 
$\beta:=d^*\omega_2$.  Since $(d^*)^2=0$, we have $d^*\alpha = d^*\beta =0$.  
Let $h(i,j)$ be the $i,j$ minor of the inverse metric tensor matrix 
$h^{ij}$.  Also, 
given distinct $j,k\in\{1,2,3\}$, let $i(j,k)$ be the number in $\{1,2,3\}$ 
that is neither $j$ nor $k$.  Also, define $\mathrm{sgn}(j,k,i(j,k))$ as $+/-1$, whichever satisfies the 
equality
\begin{equation*}
dy_1\wedge dy_2\wedge dy_3 = \mathrm{sgn}(j,k,i(j,k))dy_j\wedge dy_k\wedge dy_{i(j,k)}.
\end{equation*}
Then one can check that 
\begin{equation}
*(dy_j\wedge dy_k) = \mathrm{sgn}(j,k,i(j,k))a\cdot\sum_l (-1)^{l+i(j,k)}h(i(j,k),l)dy_l .
\end{equation} 
So we calculate
\begin{eqnarray*}
\alpha &=& d^*(\omega_1) = *d*(\omega_1)=*d(-F\cdot Ady_1)\\
&=& *(F_2 A dy_1\wedge dy_2 - F_3 A dy_3\wedge dy_1) \\
&=& (a\cdot\sum_i (-1)^{i+1}(F_2h(3,i)+ F_3h(2,i))dy_i)A.
\end{eqnarray*}
Note that in our Type B coordinates we have $h(3,1)=h(3,2)=0$ always. So since 
$F_3 = 0$ on the boundary by (\ref{bco5}), $\alpha$ satisfies CBC.
Similarly,
\begin{eqnarray*}
\beta &=& *d(G\cdot B dy_2)= *(G_1 B dy_1\wedge dy_2 - G_3 B dy_2\wedge dy_3) \\
&=& (a\cdot\sum_j(-1)^{j+1}(G_1h(3,j)-G_3h(1,j))dy_j)B.
\end{eqnarray*}
Since  $v_3(y_3)$ is constant for $y_3\in[0,d]$, we have $G_3|(0,1)\times(0,1)
\times[0,d] = 0$. This and the fact that $h(3,1)=h(3,2)=0$ show that $\beta$ 
satisfies CBC.

To calculate $[\alpha\cdot\beta]$, we first note that by Laplace expansions 
of determinants, we have
\begin{equation}
\sum_j(-1)^{i+j}h^{kj}h(i,j) = \det(h^{ij})\delta_{ik}.
\end{equation}
Indeed, if $i=k$ then the above sum is the Laplace expansion along the $k^{th}$ row 
of $h^{ij}$. If $i$ and $k$ are distinct, then the sum is a determinant of a matrix with 
a repeated row, and thus equal to $0$.  So using the above and the fact that $h^{ij}=h^{ji}$,
we have
\begin{eqnarray*}
[\alpha\cdot\beta] &=& a^2(\sum_{i,j}(-1)^{i+j}(F_2h(3,i)+F_3g(2,i))(G_1g(3,j) - 
G_3g(1,i))h^{ij})[A,B] \\
&=& a^2(\sum_{i,j}(-1)^{i+j}(F_2G_1h(3,j)h^{ji}h(3,i) 
-F_2G_3h(1,j)h^{ji}h(3,i) + \\
& & F_3G_1h(3,j)h^{ji}h(2,i) - F_3G_3h(1,j)h^{ji}h(2,i))[A,B] \\
&=& a^2\det(h^{ij})\sum_j (F_2G_1h(3,j) - F_2G_3h(1,j))\delta_{j3} + \\ 
& & (F_3G_1h(3,j) - F_3G_3h(1,j))\delta_{j2}[A,B] \\
&=& (F_2G_1h(3,3) - F_2G_3h(1,3) + F_3G_1h(3,2) - F_3G_3h(1,2))[A,B] \\
&=& (F_2G_1h(3,3) - F_3G_3h(1,2))[A,B].
\end{eqnarray*}

Since $v_3(y_3)$ is constant on $[0,d]$, we have $G_3|(0,1)\times(0,1)\times
[0,d] = 
0$.  Since $supp(F)\subset(c,d)\times(j,d)\times[0,d)$, we have 
$F_3|(0,1)\times(0,1)\times[d,1] = 0$.  Hence, $F_3G_3\equiv 0$.  So, 
continuing the above, we have
\begin{eqnarray*}
[\alpha\cdot\beta] &=& (F_2G_1h(3,3))[A,B] \\
&=& (F_2G_1 \frac{\det(h^{ij})}{h^{33}})[A,B] \\
&=& (h^{33}a^2\psi\frac{\det(h^{ij})}{h^{33}})[A,B] \\
&=& (\det(h_{ij})\det(h^{ij})\psi)[A,B] = \psi[A,B],
\end{eqnarray*}
as desired.
\end{proof}
%
%
We now extend this to a global result, and prove our main lemma.
\begin{lem}\label{curve_im}
Let $O\subset\mathbb{R}^3$ be a bounded open set. Let $f\in C^\infty
(O\otimes\k)$.  Then
\begin{equation}\label{mainderivative}
df(\nu) = -2\tau f \mbox {  on $\partial O$}
\end{equation}
if and only if
\begin{equation*}
f\in \Span\{[\alpha\cdot\beta]: \alpha,\beta\in C^\infty_c(\Lambda^1
(O\otimes\k)), d^*\alpha =d^*\beta = 0, \mbox{$\alpha,\beta$ satisfy CBC}\}.
\end{equation*}
\end{lem}
\begin{proof}
The backward direction has already been shown in Lemma \ref{bct}.  For the 
forward direction, suppose $f$ satisfies $df(\nu) = -2\tau f \mbox {  on 
$\partial O$}$.  There exists a finite cover $\{U_k\}_{k=0}^m$ of $\bar{O}$ 
that satisfies the following: $\{U_k\}_{k=1}^m$ covers the boundary and each 
$U_k$ for $k\ge 1$ is a cube in Type A coordinates, and there is a 
partition of unity $\{\lambda_k\}_{k=0}^m$ subordinate to $\{U_k\}_{k=0}^m$ so 
that $d\lambda_k(\nu)=0$ on the boundary.  

Indeed, cover $\partial O$ with finite Type A coordinate neighborhoods 
$\{U_k\}_{k=1}^m$ where $U_k=(-\delta,\delta)^2\times[0,\delta)$.  Let $W_k:=(-\delta/4,\delta/4)^2
\times[0,\delta/4)$, and $V_k:=(-\delta/2,\delta/2)^2\times[0,\delta/2)$.  Choose a smooth $\tilde{
\gamma_k}:(-\delta,\delta)^2\to\mathbb{R}$ so that $\tilde{\gamma_k}|_{[-\delta/4,\delta/4]}
\equiv 1$ and $supp(\tilde{\gamma_k})\subset [-\delta/2,\delta/2]$.  Take a smooth 
$\eta :[0,\delta)\to\mathbb{R}$ so that $\eta|_{[0,\delta/4]}\equiv 1$ and $supp(\eta)
\subset [0,\delta/2)$.  Then set $\gamma_k:U_k\to\mathbb{R}$ as
\begin{equation*}
\gamma_k(x_1,x_2,x_3)=\tilde{\gamma_k}(x_1,x_2)\eta(x_3).
\end{equation*}
Then $\gamma_k|_{W_k}\equiv 1$ and $supp(\gamma_k)\subset V_k$.  Now, take 
open sets $W_0$ and $U_0$ of $0$ so that $\bar{W_0}\subset U_0$, $\bar{U_0}
\subset O$, and $\{W_k\}_{k=0}^m$ covers $O$.  Take a smooth function 
$\gamma_0:O\to\mathbb{R}$ such that $\gamma_0|_{W_0}\equiv 1$ and $supp
{\gamma_0}\subset U_0$.  Set $\gamma:=\sum_{k=0}^m \gamma_k$, and let 
$\lambda_k:= \gamma_k/\gamma$.  Then $\{\lambda_k\}$ is a partition of unity 
with respect to $\{U_k\}_{k=0}^m$.  Also, on the support of $\gamma_k$ for 
$k>0$,
\begin{equation*}
d\gamma_k(\nu) = \frac{\partial}{\partial x_3}(\eta(x_3))\tilde{\gamma_k}
(x_1,x_2) = 0.
\end{equation*}
Hence, $d\gamma_{k}(\nu)=0$.  For $k=0$, $\gamma_k$ vanishes in a neighborhood 
of the boundary, so $d\gamma_{k}(\nu)=0$ for all $k$.  Hence, 
$d\gamma(\nu)=0$.  So, for $k>0$,
\begin{equation*}
d\lambda_k(\nu) := \frac{d\gamma_k(\nu)\gamma-\gamma_kd\gamma(\nu)}{\gamma}=0,
\end{equation*}
as we desired.

With such a partition of unity, we have 
$d(\lambda_k f)(\nu) = -2\tau\lambda_k f \mbox {  on $\partial O$}$.  So, 
by Lemmas \ref{noboundary2} and \ref{lbo} there exists 
$\{\alpha_i\}_{i=1}^{n}, \{\beta_i\}_{i=1}^{n}$ such that each 
$\alpha_i,\beta_i\in C^\infty_c(\Lambda^1(U_k\otimes\k))$, $d^*\alpha_i=
d^*\beta_i=0$, $\alpha,\beta$ satisfy CBC, and
$\lambda_k\cdot f = \sum_{i=1}^{n}[\alpha_i\cdot\beta_i]$ on $U_k$.  
Extending the $\alpha_i$'s and $\beta_i$'s by zero, we have 
$\alpha_i,\beta_i\in C^\infty_c
(\Lambda^1(O\otimes\k))$, $d^*\alpha=d^*\beta=0$, $\alpha,\beta$ satisfy CBC,
and $\lambda_k\cdot f = \sum_{i=1}^{n}[\alpha_i\cdot\beta_i]$ on $0$.  
Thus, $\lambda_k\cdot f\in
Span\{[\alpha\cdot\beta]: \alpha,\beta\in C^\infty_c(\Lambda^1(O\otimes\k)), 
d^*\alpha =d^*\beta = 0,\mbox{$\alpha,\beta$ satisfy CBC}\}$.  So, 
$f=\sum_{k=1}^m (\lambda_k\cdot f)\in
Span\{[\alpha\cdot\beta]: \alpha,\beta\in C^\infty_c(\Lambda^1(O\otimes\k)), 
d^*\alpha =d^*\beta = 0,\mbox{$\alpha,\beta$ satisfy CBC}\}$, as desired.
\end{proof}
Recasting this with our operator $T_A$, we have
\begin{cor}\label{corequiv}
Suppose $P=\bar{O}\times K\to \bar{O}$ with the flat connection $\nabla_0$ as the base connection. 
Let $g\in\Lie(\Gk)$ be smooth.  Then
\begin{equation*}
g\in \mathrm{ker}(T_0)\mbox{ if and only if } f\in\mathcal{L}_0.
\end{equation*}
\end{cor}
\begin{proof}
Set $f=\Delta g$, and apply Lemma \ref{curve_im} to $f$.
\end{proof}
\end{section}
%
%
%
%
\begin{section}{The Generation of the Smooth Gauge Algebra}
In this section we will use brackets of the image of the curvature form to get 
every smooth function in $\Lie(\Gk)$ for the special case $P=\bar{O}\times K\to \bar{O}$.  The main tool will be 
Lemma \ref{curve_im}.  The first thing we must do is see how \eqref{mainderivative} changes when 
we introduce brackets.  More specifically, note that if $g\in\mathcal{L}_0$, then 
Lemma \ref{curve_im} says that
\begin{equation}\label{mainderivative2}
d(\Delta g)(\nu)=-2\tau\Delta g.
\end{equation}
We want to know how \eqref{mainderivative2} changes if $g$ above is 
replaced by $[g_1,g_2]$, for $g_i\in \mathcal{L}_0$.  Indeed, we have
\begin{lem}\label{smooth1}
Suppose $g_1,g_2\in\mathcal{L}_0$.  Then we have
\begin{equation}
d(\Delta([g_1,g_2]))(\nu)=-2\tau\Delta[g_1,g_2] + 3[\Delta g_1, d g_2(\nu)] 
+3[d g_1(\nu), \Delta g_2].
\end{equation}
\end{lem}
\begin{proof}
First note that
\begin{eqnarray*}
\Delta([g_1,g_2]) = [\Delta g_1,g_2]+[g_1,\Delta g_2]-2[d g_1\cdot d g_2].
\end{eqnarray*}
So we have
\begin{eqnarray}\label{brack1}
d(\Delta([g_1,g_2]))(\nu) &=& d([\Delta g_1,g_2]+[g_1,\Delta g_2]-2[d g_1\cdot d g_2])(\nu)\\
&=& [d(\Delta g_1)(\nu),g_2] + [\Delta g_1,dg_2(\nu)] + [dg_1(\nu),\Delta g_2] \\
& & + [g_1,d(\Delta g_2)(\nu)] - 2d([d g_1\cdot d g_2])(\nu).
\end{eqnarray}
By Lemma \ref{curve_im}, we have
\begin{equation}\label{brack2}
d(\Delta g_i)(\nu) = -2\tau\Delta g_i .
\end{equation}
Examining the proof of Lemma \ref{bct}, we see that if $\alpha,\beta\in H^k_{con}(\kb)$ but  
are not necessarily horizontal, then we generally have
\begin{equation}
d_A([\alpha\cdot\beta])(\nu) = -2\tau[\alpha\cdot\beta] - [d_A^*\alpha,\beta(\nu)] - [\alpha(\nu), d_A^*\beta].
\end{equation}
The above yields
\begin{equation}\label{brack3}
-2d([d g_1\cdot d g_2])(\nu) = -2(-2\tau[d g_1\cdot d g_2] - [\Delta g_1, d g_2(\nu)] - [d g_1(\nu),\Delta g_2]).
\end{equation}
Plugging in \eqref{brack2} and \eqref{brack3} into \eqref{brack1}, we have
\begin{eqnarray*}
d(\Delta([g_1,g_2]))(\nu) &=& -2\tau[\Delta g_1,g_2] + [\Delta g_1,dg_2(\nu)] + [dg_1(\nu),\Delta g_2]-2\tau[g_1,\Delta g_2] \\ & &-2(-2\tau[d g_1\cdot d g_2] - [\Delta g_1, d g_2(\nu)] - [d g_1(\nu),\Delta g_2])\\
&=& -2\tau([\Delta g_1,g_2]+[g_1,\Delta g_2]-2[d g_1\cdot d g_2]) + 3[\Delta g_1, d g_2(\nu)] \\
& & + 3[d g_1(\nu),\Delta g_2] \\
&=& -2\tau\Delta([g_1,g_2]) + 3[\Delta g_1, d g_2(\nu)] + 3[d g_1(\nu),\Delta g_2],
\end{eqnarray*}
as desired.
\end{proof}
We will now show that the new term in Lemma \ref{smooth1} is actually very 
general.
\begin{lem}\label{smooth2}
Let $F$ be a smooth $\mathfrak{k}$-valued function on $\partial O$.  Then 
there exists smooth 
$g_i,h_i\in\mathcal{L}_0$ such that
\begin{equation*}
d(\Delta(\sum_i [g_i,h_i]))(\nu)+2\tau\Delta(\sum_i[g_i,h_i])= F.
\end{equation*}
\end{lem}
\begin{proof}
Since $\mathfrak{k}$ is semi-simple, there exists $A_i,B_i,C_i\in
\mathfrak{k}$ that  
\begin{equation*}
F =\sum_i f_i[[A_i,B_i],C_i]
\end{equation*}
for some real valued smooth functions $f_i$.  So, without loss of generality, 
assume that $F=f[[A,B],C]$ for some $A,B,C\in\k$.

Take any non-negative, nonzero, real-valued $\phi\in C^\infty_c(O)$.  By the Strong Minimum principle, we have $G\phi > 0$ in $O$.  Thus, we can apply Lemma 3.4 of \cite{GT} to get 
\begin{equation*}
\frac{\partial (G\phi)}{\partial \nu} < 0.
\end{equation*}
In particular, $d(G\phi)(\nu)$ never vanishes.  We set $h:= G\phi\cdot C$.  Since $\Delta h = \phi\cdot C$ has 
compact support, $h\in\mathcal{L}_0$ by Lemma \ref{curve_im}.

Let $\{U_k\}_{k=0}^{m}$ be an open cover of $O$ such that $\{U_k\}_{k=1}^{m}$ covers $\partial O$ and $U_k$ are cubes 
in Type A coordinates for $k\ge 1$.  Let $\{\lambda_k\}_{k=1}^m$ be the corresponding partition of unity for the 
cover $\{U_k\cap\partial O\}_{k=1}^m$ of the boundary.  
We set
\begin{equation*}
f_k:=\lambda_k\cdot \frac{f}{3d(G\phi)(\nu)}.
\end{equation*}
In the cube of $U_k$, suppose the $x_3$ interval is $[0,a]$.  Choose a $C^
\infty$ function $\eta:[0,a]\to[0,1]$ such that $\eta|_{[0,a/4]}\equiv 
1$ and $supp(\eta)\subseteq([0,a/2])$. We can extend $f_k$ to a function 
$\tilde{f}$ on $U_k$ by
\begin{equation*}
\tilde{f}(x_1,x_2,x_3) = f_k(x_1,x_2)\eta(x_3)\exp(-2\tau(x_1,x_2)x_3).
\end{equation*}
Note that the support of $\tilde{f}$ lies in $U_k$, so $\tilde{f}$ is a function on all of $\bar{O}$.  On $U_k$, we 
have
\begin{eqnarray*}
d\tilde{f_k}(\nu) &=& \frac{\partial}{\partial x_3}|_{x_3=0}f_k(x_1,x_2)\eta(x_3)\exp(-2\tau(x_1,x_2)x_3) \\
&=& -2\tau(x_1,x_2)f_k(x_1,x_2)\eta(x_3)\exp(-2\tau(x_1,x_2)x_3)\\
&=& -2\tau \tilde{f_k}.
\end{eqnarray*}
By Lemma \ref{curve_im}, the above shows that $G\tilde{f_k}[A,B]\in\mathcal{L}_0$.  Let $g=\sum_{k=1}^m G\tilde{f_k}[A,B]$.
We now verify that $g$ and $h$ were well-chosen.  By Lemma \ref{smooth1} and since $\Delta h|_{\partial O}\equiv 0$,
\begin{eqnarray*}
d(\Delta([g,h]))(\nu) + 2\tau\Delta[g,h] &=& 3[\Delta g, d h(\nu)] 
+3[d g(\nu), \Delta h] \\
&=& 3[\Delta g, d h(\nu)] \\
&=& 3[\sum_k f_k[A,B], dG\phi(\nu)\cdot C] \\
&=& 3[(\sum_k \lambda_k)\frac{f}{3dG\phi(\nu)}[A,B],dG\phi(\nu)C] \\
&=& f[[A,B],C]]=F,
\end{eqnarray*}
proving the lemma.
\end{proof}
We are now at the point where we can prove our main theorem.  Let $\mathcal{F}$ be the Lie algebra generated by $\mathrm{Span(Im}(\mathcal{R}_0))$.
\begin{thm}\label{mainthm}
Suppose our principal bundle is $\bar{O}\times K\to \bar{O}$, where $O\subseteq\mathbb{R}^3$ is open and bounded.  Suppose $g\in\mathrm{Lie}(\Gk)$ and is $C^\infty$.  Then $g\in\mathcal{F}$.
\end{thm}
\begin{proof}
Let $g\in\mathrm{Lie}(\Gk)\cap C^\infty$.
Recall our linear map $T_0:C^\infty(O,\mathfrak{k})\to C^\infty(\partial O,\mathfrak{k})$ by
\begin{equation*}
T_0(f) = d(\Delta f)(\nu) + 2\tau\Delta f.
\end{equation*}
Set $u:=T_0(g)$.  By Lemma \ref{smooth2}, there exists a smooth function $f\in\mathcal{F}$ such that $T_0(f)=u$.  Since $T_0$ is 
linear, we have that $T_0(g-f)=0$.  By Lemma \ref{curve_im}, we know that $g-f\in\mathrm{Span(Im}(\mathcal{R}_0))\subseteq\mathcal{F}$.  Hence, $g=f+(g-f)\in\mathcal{F}$, as we desired.
\end{proof}
The above theorem gives us our main result.
\begin{cor}
Suppose our principal bundle is $\bar{O}\times K\to \bar{O}$, where $O\subseteq\mathbb{R}^3$ is open and bounded, and 
suppose $\nabla_{A_0}=\nabla_0$.  The holonomy group $\Hko(\nabla_0)$ with base point $\nabla_0$ of the Coulomb connection of the associated bundle $\Ck\to\Ck/\Gk$ is dense in the 
connected component of the identity of $\Gk$.
\end{cor}
Before we prove this corollary, we should mention what we mean by ``holonomy group.''  We define $\Hko(\nabla_0)$ the 
the same way it would be definied in finite dimensions.  That is $g\in\Hko(\nabla_0)$ if and only if $\nabla_0\cdot g$ can be 
connected to $\nabla_0$ by a horizontal path in $\Ck$.  It has been shown that with this definition, $\Hko(\nabla_0)$ is a Banach Lie group, and the restricted holonomy group $(\Hko)^0(\nabla_0)$ is also a Banach Lie group (for the statement of this theorem, see \cite{Vas}).
\begin{proof}
This follows directly from Lemma 7.6 and Proposition 7.7 in \cite{NR}.  Specifically, Lemma 7.6 and the beginning of the proof 
of Proposition 7.7 of \cite{NR} imply that every element of $\mathcal{F}$ is the tangent vector to a curve in $(\Hko)^0(\nabla_0)$.  Then Proposition 7.7 of \cite{NR} tells us that $(\Hko)^0(\nabla_0)$ is dense in the connected component of $\Gk$ since $\mathcal{F}$ is dense in $\mathrm{Lie}(\Gk)$, completing the proof.
\end{proof}
\end{section}
\end{chapter}
\nocite{*}

\bibliography{Will-ref}
\bibliographystyle{plain}

\end{document}